\documentclass[a4paper,11pt]{amsart}

\newtheorem{theorem}{Theorem}[section]
\newtheorem{lemma}[theorem]{Lemma}

\newtheorem{proposition}[theorem]{Proposition}

\theoremstyle{definition}
\newtheorem{remark}[theorem]{Remark}
\newtheorem{definition}[theorem]{Definition}

\numberwithin{equation}{section}

\newcommand{\Ch}{\mathrm{Ch}}

\newcommand{\hphi}{\hat\varphi}
\newcommand{\hpsi}{\hat\psi}
\newcommand{\hrho}{\hat\rho}

\newcommand{\cC}{\mathcal{C}}

\newcommand{\cG}{\mathcal{G}}

\newcommand{\cL}{\mathcal{L}}

\newcommand{\cP}{\mathcal{P}}

\newcommand{\cU}{\mathcal{U}}

\newcommand{\cX}{\mathcal{X}}
\newcommand{\cY}{\mathcal{Y}}
\newcommand{\cZ}{\mathcal{Z}}

\newcommand{\bbR}{\mathbb{R}}
\newcommand{\p}{\mathbb{P}}
\newcommand{\Q}{\mathbb{Q}}

\renewcommand{\(}{\left(}
\renewcommand{\)}{\right)}
\renewcommand{\tilde}{\widetilde}

\DeclareFontFamily{U}{rcjhbltx}{}
\DeclareFontShape{U}{rcjhbltx}{m}{n}{<->rcjhbltx}{}
\DeclareSymbolFont{hebrewletters}{U}{rcjhbltx}{m}{n}

\let\aleph\relax\let\beth\relax
\let\gimel\relax\let\daleth\relax

\DeclareMathSymbol{\aleph}{\mathord}{hebrewletters}{39}
\DeclareMathSymbol{\beth}{\mathord}{hebrewletters}{98}
\DeclareMathSymbol{\gimel}{\mathord}{hebrewletters}{103}
\DeclareMathSymbol{\daleth}{\mathord}{hebrewletters}{100}

\DeclareMathSymbol{\lamed}{\mathord}{hebrewletters}{108}
\DeclareMathSymbol{\mem}{\mathord}{hebrewletters}{109}
\DeclareMathSymbol{\ayin}{\mathord}{hebrewletters}{96}
\DeclareMathSymbol{\tsadi}{\mathord}{hebrewletters}{118}
\DeclareMathSymbol{\qof}{\mathord}{hebrewletters}{113}
\DeclareMathSymbol{\shin}{\mathord}{hebrewletters}{152}

\usepackage{comment}
\usepackage{cjhebrew}
\usepackage{amsmath}
\usepackage{amsfonts,dsfont}
\usepackage[T1]{fontenc}
\usepackage{amssymb}
\usepackage{nicefrac}
\usepackage{framed}
\usepackage{color}
\usepackage{pstricks}
\usepackage{ipa}
\usepackage{textcomp}
\usepackage{url}
\usepackage{stmaryrd}
\usepackage{mathabx}
\usepackage{xcolor}
\usepackage{graphicx}
\usepackage{caption}
\usepackage{subcaption}
\usepackage[natbibapa]{apacite}
\usepackage{tikz,tikz-cd}
\usetikzlibrary{decorations.markings}
\tikzset{negated/.style={
		decoration={markings,
			mark= at position 0.5 with {
				\node[transform shape] (tempnode) {$\backslash$};
			}
		},
		postaction={decorate}
	}
}
\tikzset{negatedEq/.style={Leftrightarrow,
		decoration={markings,
			mark= at position 0.5 with {
				\node[transform shape] (tempnode) {$\backslash$};
			}
		},
		postaction={decorate}
	}
}
\usetikzlibrary{trees,matrix}

\usepackage[left=2.8cm,top=3cm,right=2.8cm,bottom=2.5cm,head=3cm,headsep=1cm, foot=3.5cm]{geometry}
\usepackage[bottom,flushmargin,hang,multiple]{footmisc}
\usepackage{setspace}
\usepackage{xr-hyper}
\usepackage{hyperref}
\usepackage{commath}
\definecolor{darkblue}{rgb}{0.1,0.1,0.9}
\definecolor{darkred}{rgb}{0.9,0.1,0.1}

\hypersetup{colorlinks, allcolors= darkblue}

\allowdisplaybreaks
\makeatletter

\newcommand{\Rmnum}[1]{\expandafter\@slowromancap\romannumeral #1@}
\makeatother

\begin{document}
	
	\title[Bounds on Choquet Risk Measures]{Bounds on Choquet Risk Measures in Finite\vspace{0.2cm}\\Product Spaces with Ambiguous Marginals\\\vspace{0.6cm}}

%	\author[Mario Ghossoub, David Saunders, and Kelvin Shuangjian Zhang]
%	{Mario Ghossoub\vspace{0.1cm}\\
%		University of Waterloo\\
%		Department of Statistics and Actuarial Science\\
%		200 University Ave.\ W.\ \\
%		Waterloo, ON, N2L 3G1 \\
%		Canada\vspace{0.1cm}\\
%		\href{mailto:mario.ghossoub@uwaterloo.ca}{mario.ghossoub@uwaterloo.ca}\vspace{0.7cm}\\
%		David Saunders\vspace{0.1cm}\\
%		University of Waterloo\\
%		Department of Statistics and Actuarial Science\\
%		200 University Ave.\ W.\ \\
%		Waterloo, ON, N2L 3G1 \\
%		Canada\vspace{0.1cm}\\\href{mailto:dsaunders@uwaterloo.ca}{dsaunders@uwaterloo.ca}\vspace{0.7cm}\\
%		Kelvin Shuangjian Zhang$^*$\vspace{0.1cm}\\
%		University of Waterloo\\
%		Department of Statistics and Actuarial Science\\
%		200 University Ave.\ W.\ \\
%		Waterloo, ON, N2L 3G1 \\
%		Canada\vspace{0.1cm}\\
%		\href{mailto:ks3zhang@uwaterloo.ca}{ks3zhang@uwaterloo.ca}\vspace{1.2cm}}
	
	% \title[Bounds on Choquet Risk Measures]{Bounds on Choquet Risk Measures in Finite\vspace{0.2cm}\\Product Spaces with Ambiguous Marginals\\\vspace{0.6cm}}

	 \author[Mario Ghossoub, David Saunders, and Kelvin Shuangjian Zhang]
	 {Mario Ghossoub\vspace{0.1cm}\\
	 	University of Waterloo\\
		 \href{mailto:mario.ghossoub@uwaterloo.ca}{mario.ghossoub@uwaterloo.ca}\vspace{0.7cm}\\
		 David Saunders\vspace{0.1cm}\\
		 University of Waterloo\\
		 \href{mailto:dsaunders@uwaterloo.ca}{dsaunders@uwaterloo.ca}\vspace{0.7cm}\\
		 Kelvin Shuangjian Zhang\vspace{0.1cm}\\
		 University of Waterloo\\
		 \href{mailto:ks3zhang@uwaterloo.ca}{ks3zhang@uwaterloo.ca}\vspace{0.7cm}\\
		 \vspace{1.2cm}}

\address{{\bf Mario Ghossoub}: University of Waterloo -- Department of Statistics and Actuarial Science -- 200 University Ave.\ W.\ -- Waterloo, ON, N2L 3G1 -- Canada}
\email{\href{mailto:mario.ghossoub@uwaterloo.ca}{mario.ghossoub@uwaterloo.ca}\vspace{0.2cm}}

\address{{\bf David Saunders}: University of Waterloo -- Department of Statistics and Actuarial Science -- 200 University Ave.\ W.\ -- Waterloo, ON, N2L 3G1 -- Canada}
	\email{\href{mailto:dsaunders@uwaterloo.ca}{dsaunders@uwaterloo.ca}\vspace{0.2cm}}

\address{{\bf Kelvin Shuangjian Zhang}: University of Waterloo -- Department of Statistics and Actuarial Science -- 200 University Ave.\ W.\ -- Waterloo, ON, N2L 3G1 -- Canada}
 \email{\href{mailto:ks3zhang@uwaterloo.ca}{ks3zhang@uwaterloo.ca}}

%	\thanks{$^*$~Corresponding author.\vspace{0.15cm}}
	% \thanks{
		% \textit{Key Words and Phrases:} Optimal Transport, Capacities, Non-Additive Measures, Risk Measures, Choquet Integral, Cooperative Games.\vspace{0.15cm} }

	\thanks{ \vspace{0.15cm}}
	
	\thanks{Mario Ghossoub and David Saunders acknowledge financial support from the Natural Sciences and Engineering Research Council of Canada in the form of Discovery Grants (NSERC Grant Nos.\ 2018-03961 and \ 2017-04220, respectively).\vspace{0.15cm}}
	
	\thanks{Declarations of interest: none}

	\maketitle
	
	%====================================================================================
	%====================================================================================
	%====================================================================================
	
	\begin{abstract}
		We investigate the problem of finding upper and lower bounds for a Choquet risk measure of a nonlinear function of two risk factors, when the marginal distributions of the risk factors are ambiguous and represented by nonadditive measures on the marginal spaces and the joint nonadditive distribution on the product space is unknown.  
		%Specifically, we assume given (marginal) capacities on these spaces, representing the ambiguous distributions of the risk factors, and we consider the problem of finding the joint capacity on the product space with these given marginals, which maximizes or minimizes the Choquet integral of a given portfolio loss function. 
		%
		%\vspace{0.4cm}
		%
		We treat this problem as a generalization of the optimal transport problem to the setting of nonadditive measures. We provide explicit characterizations of the optimal solutions for finite marginal spaces, and we investigate some of their properties. %Furthermore, we investigate the relationship between properties of the marginal capacities and those of the optimizers (and, more generally, capacities in the feasible set). In particular, we show that the minimizing capacity $\pi_{*}$ is balanced if and only if both marginal capacities are balanced, and we describe its core explicitly in that case. In contrast, in all but the most trivial cases, the maximizing capacity $\pi^{*}$ is not balanced. 
		%
		%\vspace{0.4cm}
		%
		We further discuss the connections with linear programming, showing that the optimal transport problems for capacities are linear programs, and we also characterize their duals explicitly. %Finally, we deliver a version of the Kantorovich duality.
		Finally, we investigate a series of numerical examples, including a comparison with the classical optimal transport problem, and applications to counterparty credit risk.
	\end{abstract}\vspace{0.2cm} 
	
	{\small \textit{Key Words and Phrases:} Risk management, Optimal Transport, Non-Additive Measures, Risk Measures, Cooperative Games.
		%Optimal Transport, Capacities, Non-Additive Measures, Risk Measures, Choquet Integral, Cooperative Games.
	}
	% \small	
	%  \textbf{\textit{Keywords---}} #1

	\vspace{0.4cm}
	
	%====================================================================================
	%====================================================================================
	%====================================================================================
	
	\section{Introduction}
	
	An important problem in the literature on credit risk management is that of determining bounds on the Credit Valuation Adjustment (CVA), that is, the price adjustment on a given derivatives portfolio to account for potential counterparty credit risk losses (e.g., \cite{RosenSaundersAlpha,GlassermanYang,RosenSaundersST,RosenSaundersCVA}). A portfolio's counterparty credit risk exposure depends on market risk factors, and the likelihood of a counterparty default depends on credit risk factors. Consequently, the computation of CVA requires the modelling of potential portfolio losses as functions of these two sets of dependent risk factors. There is a large literature on the required credit risk models (e.g., \cite{McNeilFreyEmbrechts} and the references therein). In practice, counterparty exposures often depend on a large number of risk factors (equity prices, interest rates, exchange rates, etc.), leading to several challenges with their measurement and management (e.g., \cite{BrigoMoriniPallavicini,Gregory}). 
	
	\vspace{0.125cm}
	
	Joint models of market and credit risk are, in general, very difficult to develop and estimate in practice.  
	%
	%and the value of the CVA depends on the dependence structure between %market risk and credit risk. 
	%
	Hence, even when the marginal distributions of the market and credit risk factors are known, there is still uncertainty about their joint distribution and about the ensuing CVA computation. %Glasserman and Yang 
	\cite{GlassermanYang} examine bounds on CVA arising from the uncertainty about the dependence structure. They formulated the problem of finding the worst-case CVA with respect to the dependence structure between the risk factors as an Optimal Transport (OT) problem. In related work, %Memartoluie et al.\ 
	\cite{MemartoluieSaundersWirjanto} considered in a formal way the problem of finding the worst-case Expected Shortfall (ES) of a nonlinear function of market risk and credit risk, given the marginal distributions of the factors, and they showed that in the case of finite sample spaces, the problem is equivalent to a linear program. 
	%Further information on both the CVA and ES bounding problems is contained in~\citet{MemartoluieThesis}.
	Recently, %Ghossoub et al.\ 
	\cite{GHS2023} extended the problem to general spaces and to spectral risk measures. They examined the problem of finding a worst-case spectral risk measure of a nonlinear function of two risk factors with known marginals, with respect to their dependence structure. They formulated the problem as a generalized OT problem and provided a strong duality theory similar to the Kantorovich duality in classical OT theory. 
	
	\vspace{0.125cm}
	
	OT is the subject of a large literature, dating back to the seminal work of %Monge 
	\cite{monge1781memoire} and %Kantorovich 
	\cite{kantorovich1942}. %Monge 
	\cite{monge1781memoire} considered the problem of minimizing the total cost (measured using the Euclidean distance between the source and the target) of moving one mass distribution to another among all volume-preserving maps. %Kantorovich 
	\cite{kantorovich1942,kantorovich1948} later relaxed this problem by expanding the feasible set to all measure couplings with given marginal distributions and developed a duality theory for the relaxed problem. Modern OT is a large and rapidly developing field (e.g., \cite{santambrogio2015optimal,VillaniOTOldAndNew}) with applications to several areas within mathematics (e.g., \cite{RachevRuschendorf,VillaniTopicsInOT}), and applied fields such as physics (e.g., \cite{guillen2017pointwise, mccann2018displacement}),  statistics (e.g., \cite{PanaretosZemel,zhang2020wasserstein}), economics (e.g., \cite{carlier2020existence,Galichon,mccann2019concavity}), finance (e.g., \cite{HenryLabordereMFHBook,eckstein2021robust}), and machine learning (e.g., \cite{peyre2019computational,TorresPereiraAmini}), for instance.
	
	\vspace{0.125cm}
	
	In the aforementioned literature, the marginal distributions of risk factors are assumed to be given and known, but their dependence structure is unknown. In particular, the marginals are (additive) probability measures. As a result, problems of bounding risk measures of loss functions can be formulated as (generalized) OT problems, with various cost functions, depending on the particular application. In many such applications, particularly related to the modelling of decision-making under ambiguity or vagueness in beliefs, a decision-maker's attitude toward, and sensitivity to ambiguity in beliefs is represented by monotone set functions that lack additivity. Such objects are called \textit{capacities} or \textit{nonadditive measures}. See, for example, the work of \cite{quiggin82,quiggin93,Schmeidler86, schmeidler1989subjective,yaari1987dual} for theoretical foundations. In particular, the seminal contribution of %Schmeidler 
	\cite{Schmeidler86, schmeidler1989subjective} axiomatized models of decision-making under ambiguity in which the decision-maker's preferences admit a representation in terms of an expected utility with respect to a nonadditive measure. Such expectations are defined through the notion of a Choquet integral with respect to a capacity.\footnote{Note that our use of the word capacity here is distinct from the usage in the literature on optimal transport with capacity constraints (e.g., \cite{korman2015optimal, korman2015dual, PennanenPerkkio2019}), where the ``capacity constraint'' imposes an upper bound on the density of the coupling.} We refer to %Denneberg 
	\cite{denneberg1994non} or %Marinacci and Montrucchio  
	\cite{montrucchiointroduction} for more about capacities and Choquet integration.
	
	\vspace{0.125cm}
	
	In this paper, we are interested in the problem of bounding a risk measure of a nonlinear function of two risk factors, but where (i) the marginal distributions of the risk factors are ambiguous, and represented by nonadditive measures on the marginal spaces; and, (ii) the objective function is a Choquet integral. As in %Glasserman and Yang 
	\cite{GlassermanYang}, we consider the case of two risk factors defined on finite spaces. We assume given (marginal) capacities on these spaces, representing the ambiguous distributions of the risk factors, and consider the problem of finding the joint capacity on the product space with these given marginals, which maximizes or minimizes the Choquet integral of a given portfolio loss function. We treat this problem as a generalization of the OT problem to the setting of nonadditive measures. We provide explicit characterizations of the optimal solutions for finite marginal spaces, and we investigate some of their properties. Additionally, we explore connections to linear programming and present a version of the Kantorovich duality.
	
	\vspace{0.125cm}
	
	The remainder of the paper is organized as follows. Section~\ref{prelim} presents definitions and background material needed for the rest of the paper. Section~\ref{section:FinProb} formulates the problem of bounding Choquet risk measures as an OT problem with nonadditive marginals. Section~\ref{section:main_result} presents a mathematical formulation of the OT problem for capacities, investigates properties of its feasible set, and gives characterizations and explicit formulas for its solution. In addition, we further study properties of the optimal capacities (in particular, non-emptiness of the core) in terms of the corresponding properties of the marginal capacities. The explicit formula for the core of the minimizer can be found in that section. Moreover, as in the case of measures, the OT problem for capacities can be  formulated as a linear program (see \cite{Torra} for a related result), and we characterize its dual in Section \ref{section:duality}. Section \ref{section:numerical_examples} presents numerical examples comparing our problem to the classical OT problem and illustrating its use in a counterparty credit risk application. Finally, Section \ref{section:conclusion} concludes.

	% \bibliographystyle{apacite}
	% %\bibliographystyle{plain}
	% \bibliography{References}
	
	% \vspace{0.4cm}
	% %========
	% \end{document} 

\vspace{0.4cm}
%====================================================================================
%====================================================================================
%====================================================================================

\section{Preliminaries}
\label{prelim}

\subsection{Capacities and Choquet Integration}
Denote by $B\(\Sigma\)$ the vector space of all bounded and $\Sigma$-measurable real-valued functions on a given measurable space $\(S, \Sigma\)$. Then $\(B\(\Sigma\), \|\cdot\|_{sup}\)$ is a Banach space \cite[IV.5.1]{Dunford}, where $\|\cdot\|_{sup}$ denotes the supnorm.
%$$\|Y\|_{sup} := \sup \{ \, \abs{Y\(s\)}: s \in S \} < +\infty, \ \forall \, Y \in B\(\Sigma\).$$

\vspace{0.125cm}

Let $ba\(\Sigma\)$ denote the linear space of all bounded finitely additive set functions on $\(S,\Sigma\)$. %, endowed with the usual vector space operations. 
When equipped with the variation norm $\|\cdot\|_{v}$, $ba\(\Sigma\)$ is a Banach space, and $\(ba\(\Sigma\), \|\cdot\|_{v}\)$ is isometrically isomorphic to the norm-dual of the Banach space $\(B\(\Sigma\),\|\cdot\|_{sup} \)$ (e.g., \cite[IV.5.1]{Dunford}) via the duality $\langle \phi, \lambda \rangle = \int \phi \ d \lambda, \ \forall \lambda \in ba\(\Sigma\), \ \forall \phi \in B\(\Sigma\)$. Denote by $ca\(\Sigma\)$ the collection of all countably additive elements of $ba\(\Sigma\)$. Then $ca\(\Sigma\)$ is a $\|\cdot\|_{v}$-closed (and hence complete) linear subspace of $ba\(\Sigma\)$. 
%Hence, $ca\(\Sigma\)$ is $\|\cdot\|_{v}$-complete, i.e.\ $\(ca\(\Sigma\), \|\cdot\|_{v}\)$ is a Banach space. 
Henceforth, a collection of probability measures will be called weak$^{*}$-compact if it is compact in the weak$^{*}$ topology $\sigma\(ba\(\Sigma\), B\(\Sigma\)\)$. 

\vspace{0.125cm}

\begin{definition}
	A capacity (nonadditive measure) on a measurable space $\(S, \Sigma\)$ is a finite set function $\gamma: \Sigma \rightarrow \left[0,\nu(S)\right]$ such that $\gamma\(\emptyset\)=0$ and $\gamma$ is monotone; that is, for any $A,B \in \Sigma$, $\gamma\(A\) \leq \gamma\(B\)$ whenever $A \subseteq B$. When $\gamma(S) = 1$, the capacity $\gamma$ is said to be normalized. %{\color{blue} Our notation for a generic capacity changes from $\nu$ in this section to $\gamma$ in the next section.}
	
	\vspace{0.125cm}
	
	The conjugate of a capacity $\gamma$ on $\(S, \Sigma\)$ is the finite set function $\bar\gamma: \Sigma \rightarrow \left[0,\nu(S)\right]$ defined by $\bar\gamma(A) := \gamma\(S\) - \gamma(A^{c})$, for all $A\in \Sigma$. Then $\bar\gamma$ is also a capacity, and if $\gamma$ is normalized then so is $\bar\gamma$.
	
	\vspace{0.125cm}
	
	A capacity $\gamma$ is called supermodular (resp.\ submodular) if 
	$$\gamma\(A \cup B\) + \gamma\(A \cap B\) \geq \ (\hbox{resp.} \leq) \ \gamma\(A\) + \gamma\(B\), \ \forall \, A,B \in \Sigma.$$
\end{definition}

\vspace{0.125cm}
The core of a capacity $\gamma$ on $\(S, \Sigma\)$, denoted by $\cC\(\gamma\)$, is the collection of all bounded finitely additive measures $\eta$ on $\(S, \Sigma\)$ such that $\eta\(A\) \geq \gamma\(A\)$, for all $A \in \Sigma$. When nonempty, $core\(\gamma\)$ is weak$^{*}$-compact and convex.  %{\color{blue} Since we have not explicitly assumed continuity of the capacities, shouldn't the core be defined as a set of finitely additive measures?}

\vspace{0.125cm}

\begin{definition}
	Let $\gamma$ be a capacity on $\(S, \Sigma\)$. The Choquet integral of $Y \in B\(\Sigma\)$ with respect to $\gamma$ is defined by
	\begin{equation*}
		\int Y \ d\gamma := \int_{0}^{+\infty} \gamma \(\{ s \in S: Y\(s\) \geq t \}\) \ dt + \int_{-\infty}^{0} \left[\gamma \(\{ s \in S: Y\(s\) \geq t \}\) - 1\right] \ dt, 
	\end{equation*}
	
	\noindent where the integrals are taken in the sense of Riemann. 
\end{definition}

\vspace{0.125cm}

% The Choquet integral with respect to a (countably additive) measure is the usual Lebesgue integral with respect to that measure \cite[p.\ 59]{MarinacciMontrucchio}. By a classical result of Schmeidler \cite{schmeidler86}, we can represent Choquet integrals with respect to a supermodular (resp.\ submodular) capacity as a lower (resp.\ upper) envelope of Lebesgue integrals:

% \begin{proposition}[Schmeidler \cite{schmeidler86}]
	% Let $\upsilon$ be a capacity on $\(S, \Sigma\)$ and let $Y \in B\(\Sigma\)$:
	% \begin{itemize}
		% \item If $\upsilon$ is supermodular then $\int Y \ d\upsilon = \min \Big\{\int Y dP: P \in core\(\upsilon_{1}\)\Big\}$;
		% \vspace{0.125cm}
		% \item If $\upsilon$ is submodular then $\int Y \ d\upsilon = \max \Big\{\int Y dP: P \in acore\(\upsilon_{2}\)\Big\}$.
		% \end{itemize}
	% \end{proposition}

% \vspace{0.125cm}

\begin{definition}
	Two functions $Y_{1},Y_{2} \in B\(\Sigma\)$ are said to be comonotonic if 
	$$\Big [ Y_{1}\(s\) - Y_{1}\(s^{\prime}\) \Big ] \Big[ Y_{2}\(s\) - Y_{2}\(s^{\prime}\) \Big] \geq 0, \hbox{ for all } s, s^{\prime} \in S.$$
\end{definition}

%For instance any $Y \in  B\(\Sigma\)$ is comonotonic with any $c \in \mathbb{R}$. Moreover, if $Y_{1},Y_{2} \in B\(\Sigma\)$, and if $Y_{2}$ is of the form $Y_{2} = I \circ Y_{1}$, for some Borel-measurable function $I$, then $Y_{2}$ is comonotonic with $Y_{1}$ if and only if the function $I$ is nondecreasing.

\vspace{0.125cm}

If $\gamma \in ca(\Sigma)$, then the Choquet integral with respect to $\gamma$ is the usual Lebesgue integral with respect to $\gamma$ (e.g., \cite[p.\ 59]{montrucchiointroduction}). Unlike the Lebesgue integral, the Choquet integral is not an additive operator on $B\(\Sigma\)$. However, the Choquet integral is additive over comonotonic functions. 

\vspace{0.125cm}

\begin{proposition}
	Let $\gamma$ be a capacity on $\(S, \Sigma\)$. 
	\begin{enumerate}
		\item If $\phi_{1}, \phi_{2} \in B\(\Sigma\)$ are comonotonic, then $\displaystyle\int \(\phi_{1} + \phi_{2}\) \ d\gamma = \int \phi_{1} \ d\gamma + \int \phi_{2} \ d\gamma$.
		\item If $\phi_{1}, \phi_{2} \in B\(\Sigma\)$ are such that $\phi_{1} \leq \phi_{2}$, then $\displaystyle\int \phi_{1} \ d\gamma \leq \int \phi_{2} \ d\gamma$.
		\item For all $\phi \in B\(\Sigma\)$ and all $c \geq 0$, then $\displaystyle\int c\phi \ d\gamma = c\int \phi \ d\gamma$.
		\item If $\gamma$ is submodular, then for any $\phi_{1}, \phi_{2} \in B\(\Sigma\)$, $\displaystyle\int \(\phi_{1} + \phi_{2}\) \ d\gamma \leq \int \phi_{1} \ d\gamma + \int \phi_{2} \ d\gamma$.
		\item If $\gamma$ is supermodular, then for any $\phi_{1}, \phi_{2} \in B\(\Sigma\)$, $\displaystyle\int \(\phi_{1} + \phi_{2}\) \ d\gamma \geq \int \phi_{1} \ d\gamma + \int \phi_{2} \ d\gamma$.
	\end{enumerate} 
\end{proposition}

\vspace{0.125cm}

\subsection{Risk Measures}
Risk measures are real-valued functionals defined on some collection of random variables on a given probability space. They are often used either as a quantification of riskiness of a given financial position, or as a way to determine adequate capital requirements (e.g., \cite{FollmerSchied}, \cite{McNeilFreyEmbrechts}, or \cite{RuschendorfMathematicalRiskAnalysis}). Formally, a risk measure is a mapping $\rho: \cX \to \bbR$, where $\cX$ is a prespecified collection of random variables on a given probablity space $\(S, \Sigma, \p\)$. Common properties of risk measures include:

\vspace{0.125cm}

\begin{enumerate}
	\renewcommand{\labelenumi}{\textbf{\theenumi}}
	\renewcommand{\theenumi}{R.\arabic{enumi}}
	\item (Monotonicity) $\rho(X) \leq \rho(Y)$, for all $X,Y \in \cX$ such that $X \leq Y$, $\p$-a.s. 
	\vspace{0.125cm}
	\item (Positive Homogeneity) $\rho(\lambda X) = \lambda \rho(X)$, for all $X \in \cX$ and all $\lambda\in\mathbb{R}_{+}$.
	\vspace{0.125cm}
	\item (Cash Invariance) $\rho(X+c) = \rho(X) + c$, for all $X \in \cX$ and $c\in\mathbb{R}$.
	\vspace{0.125cm}
	\item (Subadditivity) $\rho(X+Y) \leq \rho(X) + \rho(Y)$ for all $X,Y \in \cX$.
	\vspace{0.125cm}
	\item (Comonotonic Additivity) $\rho(X+Y) = \rho(X) + \rho(Y)$ for all $X,Y \in \cX$ that are comonotonic.
	\vspace{0.125cm}
	\item (Law Invariance) $\rho(X)=\rho(Y)$ when $X$ and $Y$ have the same distribution under $\p$.
\end{enumerate}

\vspace{0.125cm}

A coherent risk measure \cite{ArtznerDelbaenEberHeath,Delbaen2002} is a risk measure that satisfies Axioms {\bf{R.1-R.4}}, which are considered desirable for effective risk management. A practically relevant example of a coherent risk measure, frequently used in the banking and insurance industries, is the Expected Shortfall (ES), also known as the Conditional Value-at-Risk (CVaR). If $F^{\leftarrow}_{X}(t)$ is the left-continuous quantile 
of $X$, and $\alpha\in (0,1)$, then the expected shortfall of $X$ at the confidence 
level $\alpha$ is: 
\begin{equation*}
	\mbox{ES}_{\alpha}(X) = \frac{1}{1-\alpha}\int_{\alpha}^{1} F^{\leftarrow}_{X}(t)\, dt.
\end{equation*}
%ES is frequently used in quantitative risk management in the banking and insurance industries, forming the basis for the market risk charge in the Basel Framework~(\citetalias{BaselMarketRiskCharge}). For a continuous loss random variable, ES is simply the expected loss given that losses exceed a prescribed quantile. 
If the space $(S, \Sigma, \p)$ is nonatomic, then a coherent, comonotonic additive, and law-invariant risk measure admits a representation as a \textit{spectral risk measure} (e.g., \cite[Theorem 4.93]{FollmerSchied}, \cite{Kusuoka}, \cite[Proposition 8.18]{McNeilFreyEmbrechts}, \cite{Shapiro2013}), that is, as a risk measure of the form 
\begin{equation*} 
	\rho(X) = \int_{0}^{1} \mbox{ES}_{u}(X)\, d\Gamma(u),
\end{equation*}  
for some probability measure $\Gamma$ on $[0,1]$. Moreover, by a classical result on Choquet integration (e.g., \cite{Schmeidler86}), monotone and comonotonic additive risk measures admit a representation in terms of a Choquet integral of the form
$$\rho(X) = \int X \,d\gamma,$$
for some capacity $\gamma$ on $(S, \Sigma)$.
%ES at confidence level $\alpha \in (0,1)$ is the most common example of a spectral risk measure, corresponding to $\Gamma = \delta_{\alpha}$, the Dirac delta measure at $\alpha$.

\vspace{0.125cm}

The class of spectral risk measures is flexible and rich enough to encompass several of the most popular and practically relevant risk measures. Additionally, there is a tight relationship between spectral risk measures and the subclass of Choquet risk measures called Distortion Risk Measures (DRM). These are Choquet risk measures for which the capacity $\gamma$ is of the form $T \circ \p$, for some increasing function $T: [0,1] \to [0,1]$ such that $T(0) = 1 - T(1) = 0$. The function $T$ is called a distortion function, or a probability weighting function. Indeed, it can be shown (e.g., \cite{Acerbi,FollmerSchied,Kusuoka,McNeilFreyEmbrechts,Shapiro2013}) that a spectral risk measure also admits the representation
$$\rho(X) = \int_0^1 F^{\leftarrow}_{X}(t) \, \kappa(t) \, dt,$$

\noindent where $\kappa: [0,1)\to \mathbb{R}_{+}$ is a nonnegative and increasing function that satisfies $\displaystyle\int_{0}^{1}\kappa(t)\, dt = 1$. This function is called the spectral function. For instance, for the Expected Shortfall (ES) at level $\alpha$, the spectral function is given by $\(1-\alpha\)^{-1}\mathbf{1}_{[\alpha,1]}(t)$, that is, $\mbox{ES}_{\alpha}\(X\) = \(1-\alpha\)^{-1}\displaystyle\int_{\alpha}^{1} F^{\leftarrow}_X\(t\) \, dt$. Moreover, letting 
$$T(x) = 1 - \displaystyle\int_0^{1-x} \kappa(t) \,dt, \ \ \forall \, x \in [0,1],$$ 
it follows that $T$ is a distortion function, and it can be shown that $\rho$ is DRM with respect to $T \circ \p$, that is,
$$\rho(X) = \int X \, dT \circ \p.$$

\vspace{0.125cm}

\subsection{Finite State Spaces}
Suppose that $\cZ$ is a nonempty finite set, and let $\Sigma = 2^{\cZ}$ be the collection of all of its subsets. Throughout, we identify measures on any nonempty finite set $\cZ$ with vectors 
$v\in\bbR^{|\cZ|}$ through $v(A)=\sum_{i\in A}v_{i}$. Let $\gamma$ be a capacity on $\(\cZ, \Sigma\)$. 

%The core of $\gamma$ can be empty. For example, consider $\cZ=\{z_{1},z_{2}\}$, 
%$\gamma(\{z_{1}\})=\gamma(\{z_{2}\})=\tfrac{2}{3}$, and $\gamma(\cZ)=1$. 
%Any $v\in\cC(\gamma)$ would have to satisfy $v(\{z_{i}\})\geq \tfrac{2}{3}$, but $v(\cZ) = %v(\{z_{1}\}) + v(\{z_{2}\})=1$.

\vspace{0.125cm}

\begin{definition}
	The M{\"{o}}bius transform of a capacity $\gamma$ is defined as
	\begin{equation*} 
		m^{\gamma}(A) := \sum_{B\subseteq A} (-1)^{|A\setminus B|} \gamma(B).
	\end{equation*}
\end{definition}

\vspace{0.125cm}

The Choquet integral of a function $f$ with respect to the capacity $\gamma$ can be represented in terms of the M{\"{o}}bius transform as follows:

\begin{equation}\label{ChoquetIntegral}
	\begin{split}
		\gamma(f) 
		&= \sum_{A\subseteq \cX} m^{\gamma}(A) \bigwedge_{x\in A} f_{x} 
		= \sum_{A\subseteq \cX} \sum_{B\subseteq A} (-1)^{|A\setminus B|} \gamma(B) \bigwedge_{x\in A} f_{x} \\
		&= \sum_{B\subseteq \cX} \gamma(B) \left(\sum_{A\supseteq B} (-1)^{|A\setminus B|} \bigwedge_{x\in A} f_{x}\right) 
		= \sum_{B\subseteq \cX} K_{f}(B)\gamma(B), 
	\end{split}
\end{equation}
with 
\begin{equation}\label{eqn:K_f}
	K_{f}(B): =\sum_{A\supseteq B} (-1)^{|A\setminus B|} \bigwedge_{x\in A} f_{x},
\end{equation} 
where $f_{x} = f(x)$, and $\bigwedge_{x\in A}f_{x}$ represents the minimum of $f$ on $A$ (e.g., \cite[Theorem 4.95]{Grabisch}). See \cite{Grabisch} and \cite{montrucchiointroduction} for more information about the M{\"{o}}bius transform. 

\vspace{0.125cm}

\begin{definition}\label{def:outer_envelope+} 
	Let $\cZ$ be a nonempty finite set, and let $\cG\subseteq 2^{\cZ}$ be a collection of 
	subsets containing $\cZ$ and the empty set. Suppose that a function $G:\cG\to\bbR_{+}$ satisfies 
	$G(\emptyset)=0$, and $G(A)\leq G(B)$ whenever $A,B\in\cG$, $A\subseteq B$. 
	The capacity on $\cZ$ defined by
	\begin{equation*} 
		G^{*}(B) := \inf_{\substack{A\in\cG \\
				A\supseteq B}} G(A), \text{ for all } B \in 2^{\cZ},
	\end{equation*} 
	is called the outer envelope %"upper envelope" changed to "ceiling envelope"
	of $G$. The capacity defined by
	\begin{equation*} 
		G_{*}(B) := \sup_{\substack{A\in\cG \\
				A\subseteq B}} G(A), \text{ for all } B \in 2^{\cZ},
	\end{equation*} 
	is called the inner envelope %"lower envelope" changed to "floor envelope"
	of $G$. 
\end{definition} 

\vspace{0.125cm}

When it is necessary to make $\cG$ explicit in the notation, we will write $G^{*}(B) = G^{*}(B;\cG)$ for the outer envelope, and $G_{*}(B) = G_{*}(B;\cG)$ for the inner envelope. It is easy to see that $G_{*}\leq G^{*}$.\footnote{Fix $M\subseteq\cZ$, and $A,B\in\cG$ with $A\subseteq M\subseteq B$. Then $G(A)\leq G(B)$. Minimizing over $B$ containing $M$ yields $G(A)\leq G^{*}(M)$, and then maximizing over $A$ contained in $M$ gives that $G_{*}(M)\leq G^{*}(M)$.}

\vspace{0.125cm}

\begin{definition}
	Given nonempty finite sets $\cX,\cY$, we define $\cP_{\cX,\cY}$ to be the collection of 
	all subsets of $\cX\times\cY$ of the form $A\times B$ with $A\subseteq\cX$ and 
	$B\subseteq \cY$. We define $\cP^{*}_{\cX,\cY}$ to be the collection of 
	all subsets of $\cX\times\cY$ of the form $A\times B$ with $A\subseteq \cX$ 
	and $B\subseteq \cY$, and either $A=\cX$ or $B=\cY$ (or both). That is 
	$\cP^{*}_{\cX,\cY}$ is the collection of all sets either of the form 
	$\cX\times B$ with $B\subseteq \cY$ or $A\times\cY$ with $A\subseteq \cX$.
\end{definition}

\vspace{0.125cm}

Sets in product spaces and their projections will feature prominently in the optimal solutions of our optimization problems. The notation in the next definition will be convenient.

\vspace{0.125cm}

\begin{definition}
	For a set $M\subseteq \cX\times\cY$, define: 
	\begin{gather*}
		M_{\cX}:= \{x\in\cX : \exists z=(x,y)\in M\}, \quad M_{\cY}:= \{y\in\cY : \exists z=(x,y)\in M\},\\
		\tilde M_{\cX} := \{ x\in\cX : (x,y)\in M, \;\;\forall y\in\cY\}, \quad \tilde M_{\cY} := \{ y\in\cY : (x,y)\in M, \;\;\forall x\in\cX\}.
	\end{gather*} 
\end{definition}

\vspace{0.125cm}

\noindent It is easy to see that $\tilde M_{\cX} = ((M^c)_{\cX})^c$, and $\tilde M_{\cY} =((M^c)_{\cY})^c$. 

\vspace{0.125cm}

%\begin{definition} 
%Let $f:\cZ\to\bbR$, and let $\Gamma$ be a set of capacities on $\cZ$. We define 
%the upper and lower bounds of integrals of $f$ over the class $\Gamma$ to be: 
%\begin{equation} 
%\cU(f;\Gamma) = \sup_{\gamma\in\Gamma} \gamma(f)
%\ \ \hbox{and} \ \  \cL(f;\Gamma) = \inf_{\gamma\in\Gamma} \gamma(f).
%\end{equation} 
%\end{definition} 

\begin{definition} 
	Let $k \geq 2$ be an integer. A capacity $\gamma$ on $\cZ$ is called $k$-monotone if for any sets $A_{1},\ldots,A_{k} \in \cZ$, 
	\begin{equation*} 
		\gamma\left(\bigcup_{j=1}^{k} A_{j}\right) \geq \sum_{\substack{J\subseteq \{1,\ldots,k\} \\ J\ne \emptyset}}(-1)^{|J|+1} 
		\gamma\left(\bigcap_{j\in J} A_{j}\right).
	\end{equation*} 
	
	\vspace{0.125cm}
	
	The capacity is called $k$-alternating if the above inequality is reversed. A 2-monotone capacity is  supermodular, while a 2-alternating capacity is  submodular. If $\gamma$ is $k$-monotone for all $k\geq 2$, it is called totally monotone, and if it is $k$-alternating for all $k\geq 2$, it is called totally alternating.
\end{definition}

\vspace{0.4cm}
%====================================================================================
%====================================================================================
%====================================================================================

\section{Bounds on Choquet Risk Measures}\label{section:FinProb}

\subsection{Problem Formulation}
We consider the case of a portfolio whose loss depends on two risk factors defined on two finite spaces. We assume given (marginal) capacities on these spaces, representing the ambiguous distributions of the risk factors, and we consider the problem of finding the joint capacity on the product space with these given marginals that maximizes or minimizes the Choquet integral of a given portfolio loss function. 

\vspace{0.125cm}

Specifically, let $\cX$ and $\cY$ be non-empty finite sets, and let $X$ and $Y$ be random variables on $\cX$ and $\cY$, respectively. We are given a function $L: X(\cX) \times Y(\cY) \to \bbR$ representing the loss on a portfolio consisting of the risk factors $X$ and $Y$. The distributional uncertainty, or ambiguity, about the risk factors is represented by capacities $\mu$ on $\cX$ and $\nu$ on $\cY$, to be interpreted as ambiguous beliefs about the distributions of $X$ and $Y$, respectively. 

\vspace{0.125cm}

A joint distribution for $X$ and $Y$ is represented by a capacity on the product space $\cX \times \cY$, such that the projections onto $\cX$ and $\cY$ are $\mu$ and $\nu$, respectively. 

\vspace{0.125cm}

\begin{definition}\label{PiCh}
	Let $\cX$ and $\cY$ be nonempty finite sets, $\mu$ a capacity on $\cX$, $\nu$ a capacity on $\cY$, and $\pi$ a capacity on $\cX\times \cY$. 
	\vspace{0.125cm}
	\begin{enumerate}
		\item The marginal capacities of $\pi$ on 
		$\cX$ and $\cY$, respectively, are defined by
		\begin{equation*} 
			\pi_{\cX}(A) := \pi(A\times\cY) \ \ \hbox{and} \ \  \pi_{\cY}(B) := \pi(\cX\times B), \ \ \hbox{for all $A\subseteq\cX$ and $B\subseteq\cY$}.
		\end{equation*}
		
		% \vspace{0.125cm}
		%The notation defined below for the feasible set of the optimal transport problem for capacities is borrowed from \cite{gal2019kantorovich}.
		\item The set of all capacities $\pi$ on $\cX\times\cY$ such that $\pi_{\cX}=\mu$ and $\pi_{\cY}=\nu$ is denoted by $\Pi_{\Ch}(\mu,\nu)$. 
	\end{enumerate}
\end{definition}

\vspace{0.125cm}

We are interested in evaluating a risk measure $\rho\(L(X,Y)\)$ of the portfolio loss function in the case where $\rho$ is a Choquet integral of $L\(X,Y\)$ with respect to a capacity $\pi$ on $\cX \times \cY$:
$$\rho_{\pi}\(L(X,Y)\) = \int_{\cX \times \cY} L(X,Y) \, d\pi.$$

\vspace{0.125cm}

\noindent In our framework, while the capacities $\mu$ and $\nu$ are given, no information about the dependence structure (and hence the joint distribution) of the two risk factors is available. Therefore, computing a Choquet risk measure of the portfolio loss function is not possible without further information. A natural question that arises is whether we are able to establish upper and lower bounds on the value of such a risk measures with respect to the uncertrainty about the joint capacity $\pi \in \Pi_{\Ch}(\mu,\nu)$. Specifically, our problem is that of finding 
capacities that maximize or minimize the Choquet integral of $L(X,Y)$ among all capacities in $\Pi_{\Ch}(\mu,\nu)$:

{\small
	\begin{equation}\label{ProbForm}
		\cL(L;\Pi_{\Ch}(\mu,\nu)) := \inf_{\pi\in\Pi_{\Ch}(\mu,\nu)} \rho_{\pi}\(L(X,Y)\) 
		\leq \sup_{\pi\in\Pi_{\Ch}(\mu,\nu)} \rho_{\pi}\(L(X,Y)\) =: \cU(L;\Pi_{\Ch}(\mu,\nu)).
\end{equation}} 
%We note that, since $\pi(-f)\ne -\pi(f)$ in general, it is worthwhile to develop the theories for the minimum and maximum problems in parallel.

\vspace{0.125cm}

Problem \eqref{ProbForm} can be seen as a generalization of the optimal transport problem to the setting of nonadditive measures. %We develop this theory below {\color{blue} I think the last sentence should be removed.}

%We treat this problem as a generalization of the optimal transport problem to the setting of nonadditive measures. We provide explicit characterizations of the optimal solutions for finite marginal spaces, and we investigate some of their properties. Additionsally, we also explore connections to linear programming and deliver a version of the Kantorovich duality.

\vspace{0.4cm}
%====================================================================================
%====================================================================================
%====================================================================================

\section{The Optimal Transport Problem for Capacities}\label{section:main_result}
In this  section, we formulate the optimal transport problem for capacities. Once the problem is formulated, we investigate properties 
of the feasible set. Understanding the lattice structure of the feasible set leads immediately to explicit formulas for the optimizers. 
%We note that optimal transport problems for capacities have been investigated in a more general setting 
%by~\citet{gal2019kantorovich}; we will make comments on the relations between our results and theirs at several points. 

\vspace{0.125cm}

\begin{definition}
	Let $\cX$ and $\cY$ be non-empty finite sets, and let $u$ and $v$ be probability measures on $\cX$ and $\cY$, respectively. Denote by $\Pi_{a}(u, v)$ the set of measures on $\cX \times \cY$ that have the marginals $u$ on $\cX$ and $v$ on $\cY$. That is, 
	\begin{flalign*}
		\Pi_{a}(u,v) := \Big\{ \pi ~|~ \pi \text{ is a measure on } \cX \times \cY \text{ such that } \pi(A \times \cY) = u(A), \text{ for any } A \subseteq \cX,\\
		\text{ and } \pi(\cX \times B)  = v(B), 	\text{ for any } B \subseteq \cY.\Big\}
	\end{flalign*}
\end{definition}
Given a function $f$, the optimal transport minimization problem is:
\begin{flalign} \label{OT_min}
	\inf_{\pi \in \Pi_{a}(u,v)} \pi(f)   = \inf_{\pi\in\Pi_{a}(u,v)}\sum_{x \in \cX, y \in \cY} f(x,y) \, \pi(\{(x,y)\}).
\end{flalign}
Similarly, given a function $g$, the optimal transport maximization problem is:
\begin{flalign}\label{OT_max}
	\sup_{\pi \in \Pi_{a}(u,v)} \pi(g).
\end{flalign}

\vspace{0.125cm}

Both the maximization and minimization problems are linear in $\pi$. Because $\Pi_{a}(u,v)$ is convex and compact, optimal 
solutions exist, and the set of optimal solutions contains at least one extreme point of the feasible set. For instance, when $|\cX| = |\cY|$ and both $u$ and $v$ are uniform measures, by Birkhoff's Theorem there exists an optimal solution supported on 
$\bigcup_{i=1}^{|\cX|} \{(x_{i},y_{\sigma(i)})\}$, for some permutation $\sigma$.

\vspace{0.125cm}

\begin{definition} 
	Let $\cX$ and $\cY$ be nonempty finite sets and $\pi$ be a capacity on $\cX\times \cY$. 
	The marginal capacities of $\pi$ on 
	$\cX$ and $\cY$, respectively, are defined by
	\begin{equation*} 
		\pi_{\cX}(A) := \pi(A\times\cY) \ \ \hbox{and} \ \  \pi_{\cY}(B) := \pi(\cX\times B),
	\end{equation*}
	for all $A\subseteq\cX$, $B\subseteq\cY$.
\end{definition} 

% \vspace{0.125cm}

% The notation defined below for the feasible set of the optimal transport problem for capacities is borrowed from \cite{gal2019kantorovich}.
% \begin{definition} 
	% 	Let $\mu$ be a capacity on $\cX$ and $\nu$ a capacity on $\cY$. The set of 
	% 	all capacities $\pi$ on $\cX\times\cY$ such that $\pi_{\cX}=\mu$ and 
	% 	$\pi_{\cY}=\nu$ is denoted by $\Pi_{\Ch}(\mu,\nu)$. 
	% \end{definition}

\vspace{0.125cm}

In particular, for two probability measures $u$ and $v$, $\Pi_{a}(u,v) \subseteq \Pi_{\Ch}(u,v)$, where the latter is defined in Definition \ref{PiCh}. The proof of the following result is straightforward.

\vspace{0.125cm}

\begin{lemma}\label{ConjugateLemma}
	Let $\mu$ and $\nu$ be normalized capacities on $\cX$ and $\cY$, respectively. Then 
	$\pi\in\Pi_{\Ch}(\mu,\nu)$ if and only if $\bar\pi\in\Pi_{\Ch}(\bar\mu,\bar\nu)$.
\end{lemma}

\vspace{0.125cm}

Given a function 
$f:\cX\times\cY \to\bbR$,
consider the analogue of the optimal transport problem on capacity couplings, i.e.\ finding 
capacities to maximize or minimize the Choquet integral of $f$ among 
all capacities in $\Pi_{\Ch}(\mu,\nu)$:
\begin{equation*}
	\cL(f;\Pi_{\Ch}(\mu,\nu)) := \inf_{\pi\in\Pi_{\Ch}(\mu,\nu)}\pi(f)
	%\leq \gamma(f)
	\leq \sup_{\pi\in\Pi_{\Ch}(\mu,\nu)}\pi(f) =: \cU(f;\Pi_{\Ch}(\mu,\nu)).
\end{equation*} 
We note that, since $\pi(-f)\ne -\pi(f)$ in general, it is worthwhile to develop the theories for the minimum and maximum problems in parallel.

\vspace{0.125cm}
%====================================================================
%====================================================================

\subsection{The Feasible Set and Its Properties}
The first thing to observe about the feasible set is that it is nonempty.
\begin{proposition} Let $\mu$ and $\nu$ be normalized capacities on $\cX$ and $\cY$ respectively. Then 
	$\Pi_{\Ch}(\mu,\nu)\ne \emptyset$.
\end{proposition}

\vspace{0.1cm}

\begin{proof}
	Define the function $G:\cP_{\cX,\cY}\to\bbR_{+}$ by $G(A\times B):=\mu(A)\cdot\nu(B)$ for $A\times B\in \cP_{\cX,\cY}$ with 
	$A\subseteq\cX$ and $B\subseteq\cY$. It is easy to verify that both $G_{*}$ and $G^{*}$ are in 
	$\Pi_{\Ch}(\mu,\nu)$.
\end{proof}

\vspace{0.3cm}

We note that we could have used $\cP^{*}_{\cX,\cY}$ in place of $\cP_{\cX,\cY}$ in the above argument, and reached the same conclusion. Since $\Pi_{\Ch}(\mu,\nu)$ is defined by a finite system of linear equalities and inequalities, 
and $0\leq \pi(B)\leq 1$ for any set $B$, we in fact have the following result.

\vspace{0.125cm}

\begin{proposition} \label{PolyhedronProp}
	Let $\mu$ and $\nu$ be normalized capacities on $\cX$ and $\cY$. Then $\Pi_{\Ch}(\mu,\nu)$ is a compact, convex 
	polyhedron in $\bbR^{2^{|\cX|\cdot |\cY|}}$.
\end{proposition}

\vspace{0.125cm}

\begin{remark}\label{SpecialCapacitiesRemark}
	\hfil
	\begin{itemize} 
		\item A capacity $\gamma$ is called the unanimity game associated with the set $F$ if  $\gamma(G) = 1$ if $G\supseteq F$, and 
		$\gamma(G) = 0$ otherwise. If $\mu$ is the unanimity game associated with $A\subseteq \cX$, and $\nu$ is the unanimity game 
		associated with $B\subseteq \cY$, then the unanimity game $\pi$ associated with $A\times B \subseteq \cX\times\cY$ is in 
		$\Pi_{\Ch}(\mu,\nu)$.
		
		\vspace{0.125cm}
		
		\item Suppose that $\mu$ is a totally monotone capacity on $\cX$ with M{\"{o}}bius transform $m^{\mu}$, and $\nu$ is a totally monotone 
		capacity on $\cY$ with M{\"{o}}bius transform $m^{\nu}$, then $\pi$ defined to be the capacity on $\cX\times\cY$ with M{\"{o}}bius 
		transform given by
		\begin{equation*} 
			m^{\pi}(F) = \begin{cases} 
				m^{\mu}(A)\cdot m^{\nu}(B), & F = A\times B, A\subseteq \cX, B\subseteq \cY; \\
				0, & \mathrm{otherwise},
			\end{cases}
		\end{equation*} 
		is a totally monotone capacity in $\Pi_{\Ch}(\mu,\nu)$.\footnote{It should be noted that if $\mu$ and $\nu$ are capacities, this construction does not in general result in a capacity. A counterexample is given by $\cX=\cY=\{0,1\}$, $\mu=\nu$, with 
			$\mu(\emptyset)=0, \mu(\{0\})=\mu(\{1\})=0.7, \mu(\cX)=1$  (see \cite{Dyckerhoff}). 
		} 
		For further information on this construction, see \cite{Bauer2012,Destercke,ghirardato1997independence,hendon1991s,Koshevoy,WalleyFine}. Combining the above argument with Lemma~\ref{ConjugateLemma}, it is easy to see that if 
		$\mu$ and $\nu$ are totally alternating, then there exists a totally alternating capacity $\pi\in\Pi_{\Ch}(\mu,\nu)$.
		
		\vspace{0.125cm}
		
		\item A {\it possibility measure} $\gamma$ is defined as a normalized capacity such that $\gamma(A\cup B) = \max(\gamma(A),\gamma(B))$, for any sets $A$ and $B$. From this definition, it is easy to see that $\gamma(A) = \max_{z\in A} \gamma(\{z\})$ (and by normalization, there must exist $z$ such that $\gamma(\{z\})=1$). If $\mu$ and $\nu$ are possibility measures, then $\pi(A) := \max_{(x,y)\in A} \mu(\{x\})\cdot\nu(\{y\})$ defines a possibility measure in $\Pi_{\Ch}(\mu,\nu)$. The conjugate of a possibility measure is called a necessity measure (which satisfies $\gamma(A\cap B) = \min(\gamma(A),\gamma(B))$). Again, using Lemma~\ref{ConjugateLemma} one can show that if $\mu$ and $\nu$ are necessity measures, then there exists a necessity measure $\pi\in\Pi_{\Ch}(\mu,\nu)$.
	\end{itemize}
\end{remark}

\vspace{0.125cm}

A capacity is said to be balanced if its core is nonempty. The next result demonstrates that there exists a balanced $\pi\in\Pi_{\Ch}(\mu,\nu)$ 
if and only if both $\mu$ and $\nu$ are balanced.

\vspace{0.125cm}

\begin{proposition}\label{prop:core_equiv}
	Let $\mu$ and $\nu$ be normalized capacities on nonempty finite sets $\cX$ and $\cY$, respectively. Then the following are equivalent:
	\vspace{0.125cm}
	\begin{enumerate}
		\item Both $\mu$ and $\nu$ have nonempty cores (i.e., $\cC(\mu)\ne\emptyset$ and $\cC(\nu)\ne\emptyset$).
		\vspace{0.125cm}
		
		\item There exists $\pi\in\Pi_{\Ch}(\mu,\nu)$ with a nonempty core.
	\end{enumerate}
\end{proposition}

\vspace{0.1cm}

\begin{proof}
	Suppose that $u\in\cC(\mu)$ and $v\in\cC(\nu)$. Define a measure $w$ on $\cX\times\cY$ 
	by $w(\{(x,y)\}) := u(\{x\})v(\{y\})$ and additivity. Further, define 
	$G:\cP_{\cX,\cY}\to\bbR_{+}$ by $G(A\times B):=\mu(A)\cdot\nu(B)$ for 
	$A\subseteq\cX$ and $B\subseteq\cY$, and take $\pi=G_{*} \in \Pi_{\Ch}(\mu,\nu)$. It is easy to see 
	that $\pi(\cX\times\cY)=w(\cX\times\cY)$.
	Let $M\subseteq \cX\times \cY$, and consider $K=A\times B\in\cP_{\cX,\cY}$, $K\subseteq M$.  
	Then: 
	\begin{align*} 
		G(K) \, = \, \mu(A)\,\nu(B) 
		&\, \leq \, \sum_{x\in A}\sum_{y\in B} u(\{x\})v(\{y\}) 
		\, = \sum_{z=(x,y)\in K} w(\{(x,y)\}) \\
		&\leq \sum_{z=(x,y)\in M} w(\{(x,y)\}) 
		\,=\, w(M).
	\end{align*}
	This implies that $\pi(M) = G_{*}(M) \le w(M)$, for all $M\subseteq \cX\times \cY$. Therefore, $w\in\cC(\pi)$.
	
	\vspace{0.3cm}
	
	Conversely, let $\pi\in\Pi_{\Ch}(\mu,\nu)$ and $w\in\cC(\pi)$, and define for $y\in\cY$, 
	$v(\{y\}):=\sum_{x\in\cX}w(\{x,y\})$. With $B\subseteq\cY$, we have
	\begin{equation*} 
		v(B)=\sum_{y\in B} v(\{y\}) = \sum_{x\in\cX,y\in B} w(\{x,y\}) = 
		w(\cX\times B) \geq \pi(\cX\times B) = \nu(B),
	\end{equation*} 
	with equality when $B=\cY$, and therefore $v\in\cC(\nu)\ne\emptyset$. The same argument yields $\cC(\mu)\ne\emptyset$.
\end{proof}

\vspace{0.125cm}

\begin{remark}
	It should be noted that there can exist capacities $\mu$ on $\cX$ and $\nu$ on 
	$\cY$ with nonempty cores and an element $\pi\in\Pi_{\Ch}(\mu,\nu)$ with an empty core. 
	Consider $\cX=\{x_{1},x_{2}\}$, $\cY=\{y_{1},y_{2}\}$, and take $\mu$ and $\nu$ 
	to be probability measures on $\cX$ and $\cY$ respectively, 
	giving equal weight to each point. Define $\pi\in\Pi_{\Ch}(\mu,\nu)$ to give value zero to the empty set, 1 to $\cX\times\cY$, 
	$\tfrac{1}{4}$ to any subset consisting of a single point, $\tfrac{1}{2}$ to any 
	subset consisting of two points, and $\tfrac{7}{8}$ to any subset consisting of three points. Any element $w\in\cC(\pi)$ 
	would have to satisfy $w(\{(x_{1},y_{1})\})\geq \tfrac{1}{4}$, and 
	$w(\cX\times\cY\setminus\{(x_{1},y_{1})\})\geq \tfrac{7}{8}$, and thus 
	$w(\cX\times\cY) \geq \tfrac{9}{8} > 1$, contradicting 
	$w(\cX\times\cY)=\pi(\cX\times\cY)=1$.
\end{remark}

\vspace{0.125cm}
%====================================================================
%====================================================================

\subsection{Lattice Structure of the Feasible Set and Characterization of the Optimal Solutions}\label{section:characterization}
If we think of normalized capacities on $\cZ$ as functions on the collection of subsets $2^{\cZ}$, then given two capacities $\gamma$ and $\pi$, we can define, for $A\subseteq \cZ$: 
\begin{equation*} 
	(\pi\wedge \gamma)(A) := \min(\pi(A),\gamma(A)), \quad (\pi \vee \gamma)(A) := \max(\pi(A),\gamma(A)).
\end{equation*} 

\vspace{0.125cm}

\noindent With these definitions, $\pi\wedge\gamma$ and $\pi\vee\gamma$ are both capacities, and the collection of all normalized capacities is a bounded 
distributive lattice, with largest element giving value 1 to all nonempty sets, and smallest element giving value 0 to all sets except $\cZ$, which has value 1.\footnote{We note that there is another way of defining lattice operations on capacities, involving setwise maxima and minima of their 
	M{\"{o}}bius transforms. See \cite{Grabisch,montrucchiointroduction} for details.} 

\vspace{0.125cm}

Since all capacities in $\Pi_{\Ch}(\mu,\nu)$ have the same values for sets of the form $A\times\cY$, for $A\subseteq\cX$, and $\cX\times B$, for $B\subseteq\cY$, we have that $\Pi_{\Ch}(\mu,\nu)$ is a distributive sublattice. Furthermore, $\Pi_{\Ch}(\mu,\nu)$ is bounded (as a lattice) with maximum and minimum elements given by taking setwise maxima and minima:
\begin{equation*} 
	\pi^{*}(A) = \sup_{\pi\in\Pi_{\Ch}(\mu,\nu)} \pi(A) \ \ \hbox{and} \ \ 
	\quad \pi_{*}(A) = \inf_{\pi\in\Pi_{\Ch}(\mu,\nu)} \pi(A). \label{MaxMinCapacities}
\end{equation*} 

\vspace{0.125cm}

The next result follows from the definition of the Choquet integral.

\vspace{0.1cm}

\begin{theorem}
	For $f:\cX\times\cY\to\bbR$, and $\pi_{*}$ and $\pi^{*}$ described above, we have
	\begin{equation*}
		\min_{\pi\in\Pi_{\Ch}(\mu,\nu)}\pi(f) = \pi_{*}(f) \ \ \hbox{and} \ \ 
		\max_{\pi\in\Pi_{\Ch}(\mu,\nu)}\pi(f) = \pi^{*}(f).
	\end{equation*}
\end{theorem}

\vspace{0.1cm}

\begin{proof}
	We first verify that both $\pi_{*}$ and $\pi^{*}$ are indeed feasible. 
	Note that if $N=A\times\cY$ for $A\subseteq \cX$, then $\pi(N)=\mu(A)$ 
	for all $\pi\in\Pi_{\Ch}(\mu,\nu)$, and therefore $\pi^{*}(N)=\pi_{*}(N)=\mu(A)$. 
	Similarly, if $N=\cX\times B$ with $B\subseteq \cY$, then 
	$\pi^{*}(N)=\pi_{*}(N)=\nu(B)$. Furthermore, by their definitions, 
	both $\pi^{*}$ and $\pi_{*}$ are non-negative non-decreasing set functions, i.e. 
	capacities. In other words, we have that 
	$\pi^{*},\pi_{*} \in\Pi_{\Ch}(\mu,\nu)$.
	
	\vspace{0.3cm}
	
	Now, by the definition in~\eqref{MaxMinCapacities}, $\pi_{*}$ and $\pi^{*}$ achieve the set-wise 
	infimum and supremum among $\Pi_{\Ch}(\mu,\nu)$, respectively. Let $\pi\in\Pi_{\Ch}(\mu,\nu)$. Then: 
	\begin{align*}
		\pi(f) &= \int_{0}^{\infty} \pi(\{f\geq t\})\, dt  + \int_{-\infty}^{0} (\pi(\{f\geq t\})-\pi(\cZ))\, dt \\
		&= \int_{0}^{\infty} \pi(\{f\geq t\})\, dt  + \int_{-\infty}^{0} (\pi(\{f\geq t\})-1)\, dt \\
		&\geq \int_{0}^{\infty} \pi_{*}(\{f\geq t\})\, dt  + \int_{-\infty}^{0} (\pi_{*}(\{f\geq t\})-1)\, dt =\pi_{*}(f).
	\end{align*}
	The proof for $\pi^{*}$ is similar.
\end{proof}

\vspace{0.3cm}

It is possible to find explicit expressions for $\pi_{*}$ and $\pi^{*}$.

\vspace{0.1cm}

\begin{theorem} \label{thm:explicit_formula}
	For any $N\subseteq \cX\times\cY$, 
	$$\pi_{*}(N) = \max\left(\mu(\tilde N_{\cX}),\nu(\tilde N_{\cY})\right) \ \ \hbox{and} \ \ 
	\pi^{*}(N) = \min\left(\mu(N_{\cX}),\nu(N_{\cY})\right).$$
\end{theorem}

\vspace{0.1cm}

\begin{proof}
	Define $G: \cP^{*}_{\cX,\cY} \rightarrow \bbR$ by
	\begin{flalign*}
		G(M) := 
		\begin{cases}
			\mu(A), & \text{ if } M = A \times \cY;\\
			\nu(B), & \text{ if } M = \cX \times B.
		\end{cases}
	\end{flalign*}
	
	% For any 
	% set $N\subseteq \cX\times\cY$, we can then define: 
	% \begin{equation*} 
		% G^{*}(N)%: = G^{*}(N; \cP^{*}_{\cX,\cY}) 
		% = \inf_{\substack{M\in\cP^{*}_{\cX,\cY} \\ M\supseteq N}} G(M) 
		% %= \min(\mu(N_{\cX}),\nu(N_{\cY}))
		% ,
		% \end{equation*}
	% \begin{equation*}
		% G_{*}(N) %: = G_{*}(N; \cP^{*}_{\cX,\cY})
		% = \sup_{\substack{M\in\cP^{*}_{\cX,\cY}\\ M\subseteq N}} G(M)
		% %= \max(\mu(\tilde N_{\cX}),\nu(\tilde N_{\cY}))
		% .
		% \end{equation*} 
	
	\noindent Let $G*$ and $G_*$ be the outer and inner envelope of $G$ as defined in Definition \ref{def:outer_envelope+} with $\cG = \cP^{*}_{\cX,\cY}$. From the monotonicity of $\mu$ on $2^{\cX}$ (with the inclusion order), it is not hard to see that, for any $N\in\cP^{*}_{\cX,\cY}$ with $N=A\times\cY$, one has $G^{*}(N)= G_{*}(N)=\mu(A)$. Similarly, for any $N=\cX\times B$ with $B\subseteq \cY$, we have $G^{*}(N)=G_{*}(N)=\nu(B)$. By definition, $G^{*}$ and $G_{*}$ are clearly non-negative and non-decreasing, so 
	$G^{*}, G_{*} \in\Pi_{\Ch}(\mu,\nu)$.
	
	\vspace{0.3cm}
	
	For any $N\subseteq \cX\times\cY$, $\tilde N_{\cX}\times\cY \subseteq N \subseteq N_{\cX} \times \cY$ and  $\cX \times \tilde N_{\cY} \subseteq N \subseteq \cX \times N_{\cY}$. Therefore,
	$G^* (N) \le \min(\mu(N_{\cX}),\nu(N_{\cY}))$ and $G_{*}(N)  \ge \max(\mu(\tilde N_{\cX}),\nu(\tilde N_{\cY}))$. If $N\subseteq A\times \cY$, then $N_{\cX}\subseteq A$, and if $A'\times \cY\subseteq N$, then $A'\subseteq \tilde N_{\cX}$. The monotonicity of $\mu$ and $\nu$ then 
	imply that
	\begin{flalign*}
		G^* (N) = \min(\mu(N_{\cX}),\nu(N_{\cY})),\\
		G_{*}(N) = \max(\mu(\tilde N_{\cX}),\nu(\tilde N_{\cY})).
	\end{flalign*}
	
	\vspace{0.3cm}
	
	To complete the proof, we will show that $\pi_{*} =  G_{*}$ and $\pi^{*} = G^{*}$. For any $\pi\in\Pi_{\Ch}(\mu,\nu)$ and $N\subseteq \cX\times\cY$, the relation $\tilde N_{\cX}\times\cY \subseteq N \subseteq N_{\cX} \times \cY$ implies that
	\begin{equation*} 
		\mu(\tilde N_{\cX}) = \pi(\tilde N_{\cX}\times\cY) \leq 
		\pi(N)\leq 
		\pi(N_{\cX}\times\cY) = \mu(N_{\cX}),\\
	\end{equation*} 
	and $\cX \times \tilde N_{\cY} \subseteq N \subseteq \cX \times N_{\cY}$ implies that
	\begin{equation*} 
		\nu(\tilde N_{\cY}) 
		= \pi(\cX\times\tilde N_{\cY})
		\leq \pi(N)\leq \pi(\cX\times N_{\cY}) = 
		\nu(N_{\cY}).
	\end{equation*}
	Therefore, 
	\begin{equation*}
		G_{*}(N) = \max(\mu(\tilde N_{\cX}),\nu(\tilde N_{\cY})) \leq \pi(N)\leq \min(\mu(N_{\cX}),\nu(N_{\cY}))=G^{*}(N).
	\end{equation*}
	
	\vspace{0.1cm}
	
	\noindent This implies, $G_{*} \le \pi_{*}$ and  $\pi^* \le G^{*}$. The equalities hold because $G_{*}, G^{*} \in \Pi_{\Ch}(\mu,\nu)$.
\end{proof}

%The solution of this optimal transport surplus maximization problem is the 
%capacity defined by 
%$\pi^{*}(N) = \min(\mu(N_{\cX}),\nu(N_{\cY}))$, and the solution to the cost minimization problem has the explicit formula 
%$\pi_{*}(N) =  \max\left(\mu(\tilde N_{\cX}),\nu(\tilde N_{\cY})\right)$. 

\vspace{0.125cm}

\begin{remark}
	If we explicitly include the dependence of the optimizers on the marginal capacities, i.e. when given $\mu,\nu$ write $\pi_{*}(\cdot;\mu,\nu)$ 
	and $\pi^{*}(\cdot;\mu,\nu)$ for the smallest and largest elements of $\Pi_{\Ch}(\mu,\nu)$, then it is easy to show that 
	$\bar\pi_{*}(\cdot;\mu,\nu) = \pi^{*}(\cdot;\bar\mu,\bar\nu)$ and $\bar\pi^{*}(\cdot;\mu,\nu) = \pi_{*}(\cdot;\bar\mu,\bar\nu)$.
\end{remark}

\vspace{0.125cm}

\begin{remark} \hfil
	\begin{itemize} 
		\item Suppose that $\mu$ is the unanimity game associated with $A\subseteq\cX$ and $\nu$ is the unanimity game associated with 
		$B\subseteq\cY$, and $N\subseteq\cX\times\cY$. Then $\pi_{*}(N)=1$ if either $A\times\cY\subseteq N$ or $\cX\times B\subseteq N$, and 
		zero otherwise. On the other hand, $\pi^{*}(N)=1$ if for all $x_{0}\in A$ there exists $y(x_{0})\in \cY$ such that $(x_{0},y(x_{0}))\in N$ {\em and\/} 
		for all $y_{0}\in B$ there exists $x(y_{0})\in \cX$ such that $(x(y_{0}),y_{0})\in N$, and $\pi^{*}(N)=0$ otherwise.
		
		\vspace{0.125cm}
		
		\item Suppose that $\mu$ and $\nu$ are possibility measures, and define $M:\cX\times\cY \to [0,1]$ by $M(x,y) := \max(\mu(\{x\}),\nu(\{y\})$. 
		Then given $N\subseteq \cX\times\cY$, 
		\begin{equation*}
			\pi_{*}(N) = \max(\max_{x\in\tilde N_{\cX}} \mu(\{x\}),\max_{y\in\tilde N_{\cY}}\nu(\{y\})) = 
			\max_{(x,y)\in \tilde N_{\cX}\times\tilde N_{\cY}} M(x,y).
		\end{equation*} 
		Define $m:\cX\times\cY \to [0,1]$ by $m(x,y) := \min(\mu(\{x\}),\nu(\{y\}))$, then 
		\begin{equation*}
			\pi^{*}(N) = \min(\max_{x\in N_{\cX}}\mu(\{x\}),\max_{y\in N_{\cY}}\nu(\{y\})) = \max_{(x,y)\in N_{\cX}\times N_{\cY}} m(x,y).
		\end{equation*}
		When $\mu$ and $\nu$ are necessity measures, then $\pi_{*}$ and $\pi^{*}$ can be calculated using the previous remark.
	\end{itemize}
\end{remark}

\vspace{0.4cm}

Consider $f:\cX\times\cY\to\bbR$. For a fixed $x\in\cX$, define
\begin{gather*} 
	f_{y}(x) := \min\{f(x,y) : y\in\cY\} \ \ \hbox{and} \ \  
	f^{y}(x) := \max\{f(x,y) : y\in\cY\},
\end{gather*} 

\vspace{0.125cm}

\noindent with $f_{x},f^{x}:\cY\to\bbR$ defined similarly. Then

\begin{align*} 
	\widetilde{\{f\geq t\}}_{\cX} &= \{x\in\cX : (x,y)\in \{f\geq t\}\; \forall y \in\cY\} \\
	&= \{x\in\cX : \min_{y\in\cY} f(x,y)\geq t\} 
	= \{f_{y} \geq t\}. 
\end{align*} 
Similarly $\widetilde{\{f\geq t\}}_{\cY} = \{f_{x}\geq t\}$, and therefore
\begin{equation*}
	\pi_{*}(\{f\geq t\}) = \max(\mu(\{f_{y}\geq t\}),\nu(\{f_{x}\geq t\})),
\end{equation*}
and 
\begin{equation*} 
	\pi_{*}(f) = \int_{0}^{\infty} \max(\mu(\{f_{y}\geq t\}),\nu(\{f_{x}\geq t\}))\, dt 
	+ \int_{-\infty}^{0} (\max(\mu(\{f_{y}\geq t\}),\nu(\{f_{x}\geq t\}))-1)\, dt,
\end{equation*} 
using the fact that we have assumed $\mu$ and $\nu$ to be normalized.

\vspace{0.3cm}

Using a similar argument, 
\begin{align*} 
	\{f\geq t\}_{\cX} &= \{x\in\cX : \exists y\in\cY, f(x,y) \geq t\} \\
	&= \{x\in\cX : \max_{y\in\cY} f(x,y)\geq t\}
	= \{f^{y} \geq t\}, 
\end{align*} 
and $\{f \geq t\}_{\cY}=\{f^{x} \geq t\}$. Thus,
\begin{equation*} 
	\pi^{*}(\{f \geq t\}) = \min(\mu(\{f^{y}\geq t\}),\nu(\{f^{x}\geq t\})),
\end{equation*} 
and 
\begin{equation*} 
	\pi^{*}(f) = \int_{0}^{\infty} \min(\mu(\{f^{y}\geq t\}),\nu(\{f^{x}\geq t\}))\, dt 
	+ \int_{-\infty}^{0} (\min(\mu(\{f^{y}\geq t\}),\nu(\{f^{x}\geq t\}))-1)\, dt.
\end{equation*} 

\vspace{0.3cm}

To conclude, we have

\begin{align*} 
	&\cL(f;\Pi_{\Ch}(\mu,\nu)) \\
	&= \min_{\pi\in\Pi_{\Ch}(\mu,\nu)} \pi(f) 
	= \pi_{*}(f) \\
	&= \int_{0}^{\infty} \max(\mu(\{f_{y}\geq t\}),\nu(\{f_{x}\geq t\}))\, dt 
	+ \int_{-\infty}^{0} (\max(\mu(\{f_{y}\geq t\}),\nu(\{f_{x}\geq t\}))-1)\, dt \\
	&\leq \int_{0}^{\infty} \min(\mu(\{f^{y}\geq t\}),\nu(\{f^{x}\geq t\}))\, dt 
	+ \int_{-\infty}^{0} (\min(\mu(\{f^{y}\geq t\}),\nu(\{f^{x}\geq t\}))-1)\, dt \\
	&= \pi^{*}(f) 
	= \max_{\pi\in\Pi_{\Ch}(\mu,\nu)}\pi(f) \\
	&= \cU(f;\Pi(\mu,\nu)).
\end{align*}

\vspace{0.125cm}
%====================================================================
%====================================================================

\subsection{Balancedness and Cores of the Optimal Solutions}\label{section:cores}
Since $\pi_{*}(N)\leq \pi(N)\leq \pi^{*}(N)$, for all $N\subseteq\cX\times\cY$ and $\pi\in\Pi_{\Ch}(\mu,\nu)$, we immediately obtain the following result.

\vspace{0.125cm}

\begin{proposition} \label{prop:core_2} Let $\mu$ and $\nu$ be normalized capacities on $\cX$ and $\cY$, respectively. The following statements regarding the cores hold.
	\vspace{0.125cm}
	\begin{enumerate}
		\item If $\cC(\pi^{*})\ne\emptyset$, then $\cC(\pi)\ne\emptyset$ for all $\pi\in\Pi_{\Ch}(\mu,\nu)$. 
		\vspace{0.125cm}
		\item If $\cC(\pi_{*})=\emptyset$, then 
		$\cC(\pi)=\emptyset$ for all $\pi\in\Pi_{\Ch}(\mu,\nu)$. 
		\vspace{0.125cm}
		\item In particular, 
		$\cC(\pi_{*})\ne \emptyset$ iff $\cC(\mu)\ne\emptyset$ and $\cC(\nu)\ne\emptyset$.
	\end{enumerate}
\end{proposition}

\vspace{0.1cm}

\begin{proof}
	Suppose $p\in \cC(\pi^{*})$, then for any fixed $\pi \in \Pi(\mu, \nu)$ and any $N \subseteq \cX\times\cY$, one has $p(N) \ge \pi^{*}(N) \ge \pi(N)$, with both equalities hold at $N = \cX\times\cY$. Therefore, $p \in \cC(\pi)$. Using the same argument, one can show (2). Proposition \ref{prop:core_equiv} together with  (2) implies (3).
\end{proof}

\vspace{0.125cm}

However, $\cC(\pi^{*})$ is typically empty, as per the following result.

\vspace{0.125cm}

\begin{proposition} 
	Suppose that $\mu$ and $\nu$ are normalized capacities on $\cX$ and 
	$\cY$, respectively, and $|\cX|\geq 2$, $|\cY|\geq 2$. Then $\cC(\pi^{*})=\emptyset$.
\end{proposition}

\vspace{0.1cm}

\begin{proof}
	Let $\{A_{1},A_{2}\}$ and $\{B_{1},B_{2}\}$ be partitions of $\cX$ and $\cY$ 
	respectively, and define: 
	\begin{equation*} 
		N^{1}=(A_{1}\times B_{1})\cup (A_{2}\times B_{2}) \ \ \hbox{and} \ \  
		N^{2}=(A_{1}\times B_{2})\cup (A_{2}\times B_{1}).
	\end{equation*} 
	Then $N^{1}_{\cX}=N^{2}_{\cX}=\cX$, $N^{1}_{\cY}=N^{2}_{\cY} = \cY$, so that 
	for the disjoint sets $N^{1}$ and $N^{2}$, $\pi^{*}(N^{1})=\pi^{*}(N^{2})=1$.
\end{proof}

\vspace{0.125cm}

We can in fact explicitly identify $\cC(\pi_{*})$ in terms of $\cC(\mu)$ and $\cC(\nu)$.

\vspace{0.125cm}

\begin{proposition}
	Let $\mu$ and $\nu$ be normalized capacities on $\cX$ and $\cY$, respectively. Then
	\begin{equation*} 
		\cC(\pi_{*}) = \bigcup_{u\in\cC(\mu),v\in\cC(\nu)} \Pi_{a}(u,v).
	\end{equation*}
\end{proposition}

\vspace{0.1cm}

\begin{proof}
	Let $w\in\cC(\pi_{*})$, and for each fixed $x_{0}\in\cX$, $y_{0}\in\cY$
	define $u_{w}(\{x_{0}\}):=\sum_{y\in\cY}w(\{x_{0},y\})$, and $v_{w}(\{y_{0}\}):=\sum_{x\in\cX} w(\{x,y_{0}\})$. Clearly $w\in\Pi_{a}(u_{w},v_{w})$. Furthermore, for $A\subseteq \cX$, we have
	\begin{equation*}
		u_{w}(A) = w(A\times\cY) \geq \pi_{*}(A\times\cY) = \mu(A),
	\end{equation*} 
	since $\pi_{*}\in\Pi_{\Ch}(\mu,\nu)$. Thus, $u_{w}\in\cC(\mu)$, and similarly $v_{w}\in\cC(\nu)$.
	
	\bigskip
	
	Conversely, suppose that $w\in\Pi_{a}(u,v)$ with $u\in\cC(\mu)$ and $v\in\cC(\nu)$. Clearly, $w(\cX\times\cY) = u(\cX)=\mu(\cX)=1$. 
	Let $N\subseteq\cX\times\cY$, and note that $\tilde N_{\cX}\times\cY \subseteq N$ and $\cX\times\tilde N_{\cY} \subseteq N$.  Then
	\begin{align*} 
		\pi_{*}(N) = \max(\mu(\tilde N_{\cX}),\nu(\tilde N_{\cY})) \leq 
		\max(u(\tilde N_{\cX}),v(\tilde N_{\cY}))=
		\max(w(\tilde N_{\cX}\times\cY),w(\cX\times\tilde N_{\cY})) \leq w(N).
	\end{align*} 
	That is, $w \in \cC(\pi_{*})$.
\end{proof}

\vspace{0.125cm}

\begin{remark}
	By \cite[Corollary 2.23 (ii)]{Grabisch}, $\gamma$ is supermodular if and only if for every $A\subseteq B\subseteq \cX\times\cY$ and $z\notin B$, $\Delta_{z} \gamma(A)\leq \Delta_{z} \gamma(B)$, where $\Delta_{z} \gamma (A):=\gamma(A\cup\{z\})-\gamma(A)$, and $\Delta_{z} \gamma (B)$ is defined similarly. It is well-known that if $\gamma$ is supermodular, then $\cC(\gamma)\ne\emptyset$ (e.g., \cite[Theorem 3.15]{Grabisch}).
	
	\vspace{0.125cm}
	
	Let $\cX = \{x_{1},x_{2},x_{3}\}$, and $\cY=\{y_{1},y_{2},y_{3}\}$, and let $\mu$ be the additive (and therefore supermodular) capacity with 
	$\mu(\{x_{1}\})=\mu(\{x_{2}\})=0.1$, and $\mu(\{x_{3}\})=0.8$, with $\nu$ defined on $\cY$ in the same way. Define: 
	\begin{gather*} 
		A:= \{ (x_{1},y_{2}), (x_{1},y_{3})\} \ \ \text{and} \ \  
		B:= \{ (x_{1},y_{2}), (x_{1},y_{3}), (x_{2},y_{3}), (x_{3},y_{3})\},
	\end{gather*} 
	and $z:=(x_{1},y_{1})$. Note that $\tilde A_{\cX}= \emptyset$, $\tilde A_{\cY}=\emptyset$, so $\pi_{*}(A)=0$. Also, $\widetilde{(A\cup z)}_{\cX} = \{x_{1}\}$, 
	$\widetilde{(A\cup z)}_{\cY}=\emptyset$, so $\Delta_{z} \pi_{*} (A)= \pi_{*}(A\cup z) = 
	\mu(\{x_{1}\}) = 0.1$. Furthermore, $\tilde B_{\cX}=\emptyset$, 
	$\tilde B_{\cY}=\{y_{3}\}$, 
	$\widetilde{(B\cup z)}_{\cX}=\{x_{1}\}$, and $\widetilde{(B\cup z)}_{\cY}=\{y_{3}\}$, 
	so $\pi_{*}(B)=\pi_{*}(B\cup z) = \nu(\{y_{3}\})=0.8$, and $\Delta_{z}\pi_{*}(B) = 0$. 
	Thus, we conclude that while $\pi_{*}$ has a nonempty core, it is not 
	supermodular.
\end{remark}

\vspace{0.125cm}

\begin{definition}
	A capacity $\gamma$ on $\cZ$ is said to be exact if for every $S \in 2^{\cZ}\setminus\emptyset$, there exists a core element $p \in \cC(\gamma)$ such that $p(S) = \gamma(S)$.
\end{definition}

\vspace{0.125cm}

We have seen that $\cC(\pi^{*})$ is typically empty, so that $\pi^{*}$ will not be exact. In the case when $\mu$ and $\nu$ are exact, we may ask 
whether $\pi_{*}$ is exact. That is, we define the capacity $\tilde\pi\in\Pi_{\Ch}(\mu,\nu)$ by: 
\begin{equation*} 
	\tilde\pi(N) := \min\left\{p(N) : p\in \bigcup_{u\in\cC(\mu),v\in\cC(\nu)}\Pi_{a}(u,v)\right\}, \text{ for any } N\subseteq\cX\times\cY,
\end{equation*} 

\vspace{0.125cm}

\noindent and we ask whether $\pi_{*} = \tilde \pi$.

\vspace{0.125cm}

\begin{remark}
	In general $\tilde\pi$ as defined above need not be either submodular or supermodular. To see this, consider the case 
	$\cX = \cY = \{1,2,\ldots,n\}$ for some $n\geq 3$, with $\mu$ and $\nu$ being uniform probability measures, and 
	let $\pi'$ be the conjugate of $\tilde \pi$.\footnote{We prefer to avoid the cumbersome notation $\bar{\tilde\pi}$.} Then 
	\begin{equation*} 
		\pi'(A) = 1-\tilde\pi(A^{c}) = 1-\min_{p\in\Pi_{a}(\mu,\nu)} p(A^{c}) = \max_{p\in\Pi_{a}(\mu, \nu)}p(A).
	\end{equation*} 
	By Birkhoff's Theorem, the optimum $\tilde\pi(A)$ (and similarly $\pi'(A)$) is achieved by measures that put mass $\frac{1}{n}$ on 
	points $\{x_{i},y_{\sigma(i)}\}$ for some permutation $\sigma$.
	Consider $A_{1}=\{(1,1)\}$, $z=(n,n)$ and $B_{1} = \cX\times\cY\setminus\{z\}$. Then it is easy to see 
	that  $\Delta_{z}\pi'(A_{1}) = \tfrac{2}{n}-\tfrac{1}{n} = \tfrac{1}{n}$, while $\Delta_{z}\pi'(B_{1})=1-1=0$. Thus 
	$\Delta_{z}\pi'(A_{1}) > \Delta_{z}\pi'(B_{1})$, and $A_{1}\subseteq B_{1}$, so $\pi'$ is not supermodular (and therefore 
	$\tilde\pi$ is not submodular, see \cite[Theorem 2.20]{Grabisch}). On the other hand, consider $A_{2}=\{(1,1)\}$, 
	$B_{2}=\{(1,1),(2,1)\}$ and $z=(1,2)$. Then $\Delta_{z}\pi'(A_{2})=0$, and $\Delta_{z}\pi'(B_{2})=\frac{1}{n}$. 
	We therefore have that $(B_{2}\cup\{z\})^{c}\subseteq (A_{2}\cup\{z\})^{c}$, and 
	$\Delta_{z}\tilde\pi((B_{2}\cup \{z\})^{c}) = \Delta_{z}\pi'(B_{2}) > \Delta_{z}\pi'(A_{2}) = \Delta_{z}\tilde\pi((A_{2}\cup\{z\})^{c})$ (e.g., \cite[Theorem 2.16]{Grabisch}). Thus $\tilde\pi$ is not supermodular (and $\pi'$ is not submodular).
\end{remark}

\vspace{0.125cm}

\begin{remark}\label{rmk:non_exact}
	Let $n\geq 2$, $\cX=\{1,\ldots,n\}$ and $\cY=\cX$, and 
	take $\mu$ and $\nu$ to be two probability measures on $\cX$ that are 
	not equal. Then $\cC(\mu)=\{\mu\}$, and $\cC(\nu) = \{\nu\}$, so 
	that $\cC(\pi_{*})=\Pi_{a}(\mu,\nu)$. Notice that any 
	element of $\Pi_{a}(\mu,\nu)$ is also in $\Pi_{\Ch}(\mu,\nu)$. $\Pi_{a}(\mu,\nu)$ 
	is compact, and for any fixed $B$, $p(B) =  
	\sum_{\{x,y\}\in B} p(\{x,y\})$ is a continuous function on 
	$\Pi_{a}(\mu,\nu)$ and therefore its minimum is attained. 
	% Suppose 
	% that $\pi_{*}$ was exact, and consider the set 
	% $D = \{(1,1),(2,2),\ldots,(n,n)\}$ and $M=D^{c}$. We have that 
	% $\tilde M_{\cX}=\tilde M_{\cY} = \emptyset$, and therefore 
	% $\pi_{*}(M)=0$. Because $\pi_{*}$ has been assumed to be exact, 
	Consider the set 
	$D = \{(1,1),(2,2),\ldots,(n,n)\}$ and $M=D^{c}$. We have that 
	$\tilde M_{\cX}=\tilde M_{\cY} = \emptyset$, and therefore 
	$\pi_{*}(M)=0$.
	Suppose that $\pi_{*}$ was exact.
	Then there 
	is a $\pi\in\Pi_{a}(\mu,\nu)$ such that $\pi(M)=0$. But then 
	$\pi$ is concentrated on the diagonal $D$, contradicting the fact that 
	$\mu\ne\nu$. This implies that $\pi_{*}$ is not exact. 
\end{remark}

%\begin{remark}\label{rmk:counter-example}
%The above remark seems to present a contradiction to the 
%duality result claimed by~\citet{gal2019kantorovich}.
%%\marginpar{SJ: in \cite{gal2019kantorovich}, $c$ is assumed to be continuous.}
%There is a 
%supermodular product capacity; for any nonnegative cost function $c$, 
%the minimum of $\gamma(c)$ on $\Pi(\mu,\nu)$ is $\pi_{*}(c)$; there are 
%$\pi_{*}(c)\ne \min_{\gamma\in\cC(\pi_{*})} \gamma(c)$; but by 
%optimal transport duality, the dual of 
%$\max_{\varphi \oplus \psi\leq c} \mu(\varphi) + \nu(\psi) = 
%\min_{\gamma\in\Pi_{a}(\mu,\nu)} \gamma(c)=\min_{\gamma\in\cC(\pi_{*})} 
%\gamma(c)$, 
%which is the Choquet integral of $c$ with respect to the lower  
%envelope of $\pi_{*}$. %({\bf TODO:} Notational conflict with $\citet{Grabisch}$).
%\end{remark}

\vspace{0.4cm}
%====================================================================================
%====================================================================================
%====================================================================================

\section{Linear Programming and the Kantorovich Duality for Capacities}\label{section:duality}

In this section, we formulate the optimal transport problem for capacities as a linear program, and we present its dual. Recall that the Choquet integral of $f$ with respect to a capacity $\gamma$ on $\cZ$ can be written as
\begin{align*} 
	\int f \, d\gamma 
	% &= \sum_{A\subseteq \cZ} m^{\gamma}(A) \bigwedge_{x\in A} f_{x} \nonumber \\
	% &= \sum_{A\subseteq \cZ} \sum_{B\subseteq A} (-1)^{|A\setminus B|} \gamma(B) \bigwedge_{x\in A} f_{x} \\
	% &= \sum_{B\subseteq \cZ} \gamma(B) \left(\sum_{B\subseteq A} (-1)^{|A\setminus B|} \bigwedge_{x\in A} f_{x}\right) \\
	&= \sum_{B\subseteq \cZ} K_{f}(B)\gamma(B),
\end{align*} 
where 
\begin{equation*} 
	K_{f}(B) =\sum_{A\supseteq B} (-1)^{|A\setminus B|} \bigwedge_{x\in A} f_{x}.
\end{equation*} 

\vspace{0.125cm}

\noindent While this expression is not linear in $f$, it is linear in $\gamma$, and since the constraints defining $\Pi_{\Ch}(\mu,\nu)$ 
are all linear (see Proposition~\ref{PolyhedronProp}), the problem of minimizing $\pi(f)$ over all $\pi\in\Pi_{\Ch}(\mu,\nu)$ becomes a linear program: 

\begin{gather} 
	\min_{\pi} \sum_{B\subseteq \cX\times\cY} K_{c}(B) \pi(B) \label{eqn:primal_ini},  ~~\text{ subject to}\\
	\begin{aligned}
		\pi(G\times\cY) &= \mu(G), &\quad &\emptyset \ne G\subseteq \cX; \\
		\pi(\cX\times F) &= \nu(F), &\quad &\emptyset \ne F\subseteq \cY;  \\
		\pi(A\cup w)  &\geq \pi(A), &\quad &A\subset \cX\times \cY, w =\{(x,y)\} \notin A;  \\
		\pi(\emptyset) &= 0, &
		%\\{\red \pi & \red \geq 0 }& 
		\label{eqn:discrete_primal}
	\end{aligned}
\end{gather} 
\vspace{0.125cm}

\noindent (e.g., \cite[pp.\ 81-82]{Grabisch}). Recall that a subset $B$ of $\cX\times\cY$ is in $\cP^{*}_{\cX,\cY}$ %$\cP$
if $B = G\times\cY$ for some 
$G\subseteq \cX$ or 
$B = \cX \times F$ for some $F\in \cY$.

\vspace{0.125cm}

%\begin{theorem}\label{thm:duality_min_1}
The dual of the above linear program is given by

\begin{gather} 
	\max_{\hat\varphi,\hpsi,\hrho} \sum_{G\subseteq \cX} \hat\varphi(G)\mu(G) + \sum_{F\subseteq \cY} \hpsi(F)\nu(F) \label{eqn:dual_ini}, ~~\text{ subject to}
\end{gather}
\begin{gather}
	\begin{aligned}
		\hat\varphi(G) - \sum_{w\notin G\times \cY} \hrho(G\times\cY,w) + 
		\sum_{w\in G\times\cY} \hrho((G\times\cY)\setminus \{w\},w) & = K_{c}(G\times\cY), &\quad &\emptyset \ne G\subsetneqq \cX; %\nonumber 
		\\
		\hpsi(F) - \sum_{w\notin \cX\times F} \hrho(\cX\times F,w) + 
		\sum_{w\in \cX\times F} \hrho((\cX\times F)\setminus\{w\},w) & = K_{c}(\cX\times F),  &\quad & \emptyset \ne F\subsetneqq \cY; %\nonumber 
		\\
		\hat\varphi(\cX) + \hpsi(\cY) +\sum_{w} \hrho((\cX\times\cY)\setminus\{w\},w) & = K_{c}(\cX\times\cY); %\nonumber
		\\
		-\sum_{w\notin B} \hrho(B,w) + \sum_{w\in B} \hrho(B\setminus\{ w\},w) & = K_{c}(B),  &\quad & B\notin \cP^{*}_{\cX,\cY}; %\cP 
		%\nonumber
		\\
		\hrho & \geq 0. \label{eqn:discrete_dual}
	\end{aligned}
\end{gather}
%\end{theorem}

\vspace{0.125cm}

Let $(\hphi_{*}, \hpsi_{*}, \hrho_{*})$ be an optimal solution to (\ref{eqn:dual_ini} - \ref{eqn:discrete_dual}). Then complementary slackness implies that, for any $(A, w) \in \{ (A, w) \in 2^{\cX \times \cY} \times (\cX \times \cY): w \not\in A\}$, 
\begin{flalign}
	\hrho_{*} (A, w) \, \left(\pi_{*}(A \cup w) - \pi_{*}(A)\right) = 0. 
\end{flalign}
%\end{corollary}

\vspace{0.125cm}

\begin{remark}\label{RemDualMax}
	The dual of the maximization problem 
	\begin{flalign}\label{MaxChoquetIntegral}
		\max_{\pi \in \Pi_{\Ch}(\mu, \nu)} \pi(c)
	\end{flalign}
	is given by
	\begin{gather} 
		\min_{\hat\varphi,\hpsi,\hrho} \sum_{G\subseteq \cX} \hat\varphi(G)\mu(G) + \sum_{F\subseteq \cY} \hpsi(F)\nu(F) \label{eqn:dual_ini2}, ~~\text{ subject to}
	\end{gather}
	\begin{gather}
		\begin{aligned}
			\hat\varphi(G) - \sum_{w\notin G\times \cY} \hrho(G\times\cY,w) + 
			\sum_{w\in G\times\cY} \hrho((G\times\cY)\setminus \{w\},w) & = K_{c}(G\times\cY),  &\quad & \emptyset \ne G\subsetneqq \cX; %\nonumber 
			\\
			\hpsi(F) - \sum_{w\notin \cX\times F} \hrho(\cX\times F,w) + 
			\sum_{w\in \cX\times F} \hrho((\cX\times F)\setminus\{w\},w) & = K_{c}(\cX\times F),  &\quad & \emptyset \ne F\subsetneqq \cY; %\nonumber 
			\\
			\hat\varphi(\cX) + \hpsi(\cY) +\sum_{w} \hrho((\cX\times\cY)\setminus\{w\},w) & = K_{c}(\cX\times\cY); %\nonumber
			\\
			-\sum_{w\notin B} \hrho(B,w) + \sum_{w\in B} \hrho(B\setminus\{ w\},w) & = K_{c}(B),  &\quad & B\notin \cP^{*}_{\cX,\cY}; %\cP 
			%\nonumber
			\\
			\hrho & \leq 0. \label{eqn:discrete_dual2}
		\end{aligned}
	\end{gather}
	Suppose that $(\hphi^{*}, \hpsi^{*}, \hrho^{*})$ is an optimal solution to (\ref{eqn:dual_ini2} - \ref{eqn:discrete_dual2}). Then by complementary slackness, for any $(A, w) \in \{ (A, w) \in 2^{\cX \times \cY} \times (\cX \times \cY): w \not\in A\}$, 
	\begin{flalign*}
		\hrho^{*} (A, w) \, \left(\pi^{*}(A \cup w) - \pi^{*}(A)\right) = 0. 
	\end{flalign*}
\end{remark}

\vspace{0.125cm}

%\begin{theorem}\label{thm:duality_min_2}
\begin{remark}
	The dual of the minimization Optimal Transport problem is equivalent to the problem
	\begin{flalign} 
		\max_{L_{\varphi},L_{\psi}, \hrho} 
		\sum_{G\subseteq \cX} m^{\mu} (G) L_{\varphi}(G) 
		+ \sum_{F\subseteq \cY} m^{\nu}(F)  L_{\psi}(F),~~\text{ subject to}
		\label{eqn:dual_2}
	\end{flalign}
	\begin{flalign}
		\begin{aligned}
			L_{\varphi}(A_{\cX}) + L_{\psi}(A_{\cY}) + \sum_{D \supseteq A} \sum_{w \in A} \hrho (D\setminus\{w\}, w) &= \bigwedge_{(x,y)\in A} c(x,y), \quad \emptyset \ne A \subseteq \cX\times \cY;\\
			\hrho & \geq 0. \label{eqn:dual_2_constraint}
		\end{aligned}
	\end{flalign} 
	%\end{theorem}
	
	%\begin{proof}
	\noindent To see this, we will show, by the following change of variables\footnote{This corresponds to the situation derived from a set function $\xi_\varphi$, where $\hat\varphi = m^{\xi_{\varphi}}$, and $L_{\varphi} = \check m^{\xi_{\varphi}}$, the co-M{\"{o}}bius transform, with similar 
		conventions for $\psi$, see \cite[Table A.2, p.\  440]{Grabisch}.}
	\begin{flalign*}
		\begin{aligned} 
			\hat\varphi(G) := \sum_{B\supseteq G} (-1)^{|B\setminus G|}L_{\varphi}(B); \hspace{2cm}
			\hpsi(F) := \sum_{A\supseteq F} (-1)^{|A\setminus F|} L_{\psi}(A),
		\end{aligned}
	\end{flalign*}
	that the objectives are equal and that the constraints can be derived from each other.

	\vspace{0.125cm}
	
	First, the objective function becomes
	\begin{flalign*}
		& \sum_{G\subseteq \cX} \hat\varphi(G)\mu(G) + \sum_{F\subseteq \cY} \hpsi(F)\nu(F) \\
		= & \sum_{G\subseteq \cX}  \sum_{B\supseteq G} (-1)^{|B\setminus G|}L_{\varphi}(B) \mu(G) + \sum_{F\subseteq \cY} \sum_{A\supseteq F} (-1)^{|A\setminus F|} L_{\psi}(A) \nu(F)\\
		= & \sum_{B\subseteq \cX}  \left( \sum_{G\subseteq B} (-1)^{|B\setminus G|} \mu(G) \right) L_{\varphi}(B) + \sum_{A\subseteq \cY} \left( \sum_{F\subseteq A}  (-1)^{|A\setminus F|}  \nu(F)\right) L_{\psi}(A)\\
		= & \sum_{B\subseteq \cX}  m^{\mu}(B) L_{\varphi}(B) + \sum_{A\subseteq \cY} m^{\nu}(A)  L_{\psi}(A).
	\end{flalign*}
	
	\vspace{0.125cm}
	
	\noindent To see that the constraints are equivalent, notice that the above transformation can be inverted as
	\begin{flalign*} 
		\begin{aligned}
			L_{\varphi}(G) = \sum_{G'\supseteq G} \hat\varphi(G'); \hspace{2cm}
			L_{\psi}(F) = \sum_{F'\supseteq F} \hpsi(F').
		\end{aligned}
	\end{flalign*}

	\noindent Furthermore, for any $B\subseteq \cX\times\cY$, recall 
	\begin{equation*} 
		K_{c}(B) = \sum_{A\supseteq B} (-1)^{|A\setminus B|} \bigwedge_{(x,y)\in A} c(x,y).
	\end{equation*} 
	Using the same inversion formula, we obtain
	\begin{equation*} %\label{eqn:McB}
		%M_{c}(A):=
		\bigwedge_{(x,y)\in A} c(x,y) = \sum_{B\supseteq A} K_{c}(B).
	\end{equation*} 
	
	\vspace{0.125cm}
	
	For any non-empty set $A \subseteq \cX \times \cY$, sum all constraints with a right-hand side involving $K_{c}(B)$ with $B \supseteq A$. 
	The right-hand side term of \eqref{eqn:discrete_dual} becomes 
	\begin{equation*} 
		\sum_{B\supseteq A } K_{c}(B) = \bigwedge_{(x,y)\in A} c(x,y).
	\end{equation*} 
	
	\vspace{0.125cm}
	
	The sum of terms on the left-hand side of \eqref{eqn:discrete_dual} will yield a sum involving $\hat\varphi$, which is 
	\begin{equation*} 
		\sum_{G'\supseteq A_{\cX}} \hat\varphi(G') = L_{\varphi}(A_{\cX}),
	\end{equation*} 
	and a sum involving $\hpsi$, which is
	\begin{equation*} 
		\sum_{F'\supseteq A_{\cY}} \hpsi(F') = L_{\psi}(A_{\cY}).
	\end{equation*} 
	
	\vspace{0.125cm}
	
	Lastly, %we obtain several terms involving $\hrho$ left in \eqref{eqn:discrete_dual}. 
	denoting the sum of all terms involving $\hrho$ in \eqref{eqn:discrete_dual} by $S$, we obatin
	\begin{equation*} 
		S = J_{1} + J_{2} + J_{3} + J_{4} + J_{5} + J_{6} + J_{7}, 
	\end{equation*}
	where
	\begin{align*} 
		J_{1} &:= -\sum_{\substack{G'\supseteq A_{\cX} \\ G'\ne \cX}} \sum_{\substack{x\notin G' \\ y \in \cY}} \hrho(G'\times\cY,(x,y)); 
		& J_{2} &:= \sum_{\substack{G'\supseteq A_{\cX} \\ G' \ne \cX}} \sum_{\substack{x\in G'\\ y \in \cY}} 
		\hrho(G'\times\cY \setminus \{(x,y)\},(x,y)); \\
		J_{3} &:= -\sum_{\substack{F'\supseteq A_{\cY} \\ F'\ne \cY}} \sum_{\substack{y\notin F' \\ x \in \cX}} \hrho(\cX\times F',(x,y)); &
		J_{4} &:= \sum_{\substack{F'\supseteq A_{\cY} \\ F' \ne \cY}} \sum_{\substack{y\in F'\\ x \in \cX}} 
		\hrho(\cX\times F' \setminus \{(x,y)\},(x,y)); \\
		J_{5} &:= \sum_{x\in\cX,y\in\cY} \hrho(\cX\times\cY\setminus \{(x,y)\},(x,y)); \\
		J_{6} &:= -\sum_{\substack{  B \supseteq A  \\
				B \notin \cP^{*}_{\cX,\cY}
				%\cP 
				\\ (x,y)\notin B}} 
		\hrho(B,(x,y)); &
		J_{7} &:= \sum_{\substack{ B \supseteq A  \\
				B \notin \cP^{*}_{\cX,\cY}
				\\ (x,y)\in B}} 
		\hrho(B\setminus\{(x,y)\},(x,y)).
	\end{align*}
	
	\vspace{0.125cm}
	
	\noindent By summing the above terms, we obtain 
	\begin{flalign}\label{eqn:hat_rho}
		\begin{aligned}
			S &= -\sum_{B \supseteq A} \sum_{w \notin B} \hrho (B, w) + \sum_{D \supseteq A } \sum_{w \in D} \hrho (D\setminus\{w\}, w)\\
			&= -\sum_{B \supseteq A } \sum_{w \notin B} \hrho (B, w) + \sum_{D \supseteq A} \left[\sum_{w \in A} \hrho (D\setminus\{w\}, w) + \sum_{w \in D\setminus A} \hrho (D\setminus\{w\}, w)\right]\\
			& = \sum_{D \supseteq A} \sum_{w \in A} \hrho (D\setminus\{w\}, w). 
		\end{aligned}
	\end{flalign}
	
	\vspace{0.125cm}
	
	\noindent The last equality comes from the observation that there exists an one-to-one mapping between  $\{(B, w): A \subset B, w\notin B \}$ and $\{(D, w): A \subset D, w \in D\setminus A\}$ by the map $D := B \cup \{w\}$, and thus the first and third terms in the second line of \eqref{eqn:hat_rho} cancel out. Therefore, one can derive the equations in \eqref{eqn:dual_2_constraint} from those in \eqref{eqn:discrete_dual}. Similarly, one can also prove the opposite direction by using the above change of variables. \qedhere
	%\end{proof}
\end{remark}

\vspace{0.125cm}

\begin{remark}
	By a similar argument, the dual of the maximization Optimal Transport problem is equivalent to
	\begin{flalign} 
		\min_{L_{\varphi},L_{\psi}, \hrho} 
		\sum_{G\subseteq \cX} m^{\mu} (G) L_{\varphi}(G) 
		+ \sum_{F\subseteq \cY} m^{\nu}(F)  L_{\psi}(F), ~~\text{ subject to}
	\end{flalign}
	\begin{flalign}
		\begin{aligned}
			L_{\varphi}(A_{\cX}) + L_{\psi}(A_{\cY}) + \sum_{D \supseteq A} \sum_{w \in A} \hrho (D\setminus\{w\}, w) &= \bigwedge_{(x,y)\in A} c(x,y), \quad \emptyset \ne A \subseteq \cX\times \cY;\\
			\hrho & \leq 0. 
		\end{aligned}
	\end{flalign} 
\end{remark}

\vspace{0.4cm}
%====================================================================================
%====================================================================================
%====================================================================================

\section{Numerical Examples}
\label{section:numerical_examples}

\vspace{0.125cm}

\subsection{A Comparison with the Optimal Transport Problem for Additive Measures}\label{subsec:numerical_1}

%\vspace{0.125cm}
In this section, we compare the optimal transport problem for capacities with the classical optimal transport problem (for measures) via numerical simulations. Assume that $\cX$ and $\cY$ are two finite subsets of $\mathbb{R}$ with $|\cX|=30$ and $|\cY| = 20$, and $\mu$ and $\nu$ are probability measures on $\cX$ and $\cY$, respectively. Given the quadratic function $c(x,y) = (x-y)^2$ on $\cX \times \cY$, the classical optimal transport problem is to find 
\begin{equation*}
	\min_{\pi \in \Pi_a(\mu, \nu)} \int_{\cX \times \cY} c(x,y)\, d \pi(x,y);
\end{equation*}
while the optimal transport minimization problem for capacities seeks
\begin{equation*}
	\min_{\gamma \in \Pi_{\Ch}(\mu, \nu)} \gamma(c),
\end{equation*}
where $\gamma(c)$ represents the Choquet integral of $c$ with respect to $\gamma$. The latter problem will have a lower minimum since its feasible set is larger. 

\vspace{0.125cm}

%Now, introduce the method of calculating the optimal transport problem:
We use the Python package AMPL to solve the linear program for the classical optimal transport minimization. However, the linear program for the optimal transport for capacities is quite large when both sets have cardinality greater than 5. For example, when $|\cX|=|\cY| = 5$, the number of variables in the linear program is $33,\!554,\!432$, and the number of constraints is $419,\!430,\!437$. These numbers will become astronomical if $|\cX|$ and $|\cY|$ exceed 20. For the case when $|\cX|=30$ and $|\cY| = 20$, the number of variables for the linear program is larger than $10^{180}$, and the number of constraints is larger than $10^{183}$; while the number of variables for the classical optimal transport problem is $600$ and the number of constraints is only $50$. Therefore, in this case, solving the classical optimal transport problem using linear programming methods is still fast, but the linear program for capacities cannot be solved using numerical methods. However, using the explicit solution provided in Theorem \ref{thm:explicit_formula}, the minimum can be computed in a few seconds even when $|\cX| = |\cY| = 100$. 

% \begin{figure}[hbt!]
	% 	\centering
	% \subfloat[plot with the line y=x]{%
		%   \includegraphics[clip,width=.2\columnwidth]{images/test1.png}%
		%   \label{}
		% }
	
	% \subfloat[]{%
		%   \includegraphics[clip,width=.2\columnwidth]{images/test1_loglog.png}%
		%   \label{fig:test1b}
		% }

	% %	\caption{Comparison of the optimal values of the optimal transport problem for capacities with those of the optimal transport problem for additive measures, when both marginals are additive measures.}	
	
	%  %\label{fig:test1}
	% \end{figure}

\begin{figure}[hbt!]
	\centering
	\begin{subfigure}[b]{0.47\textwidth}
		\centering
		\includegraphics[width=\textwidth]{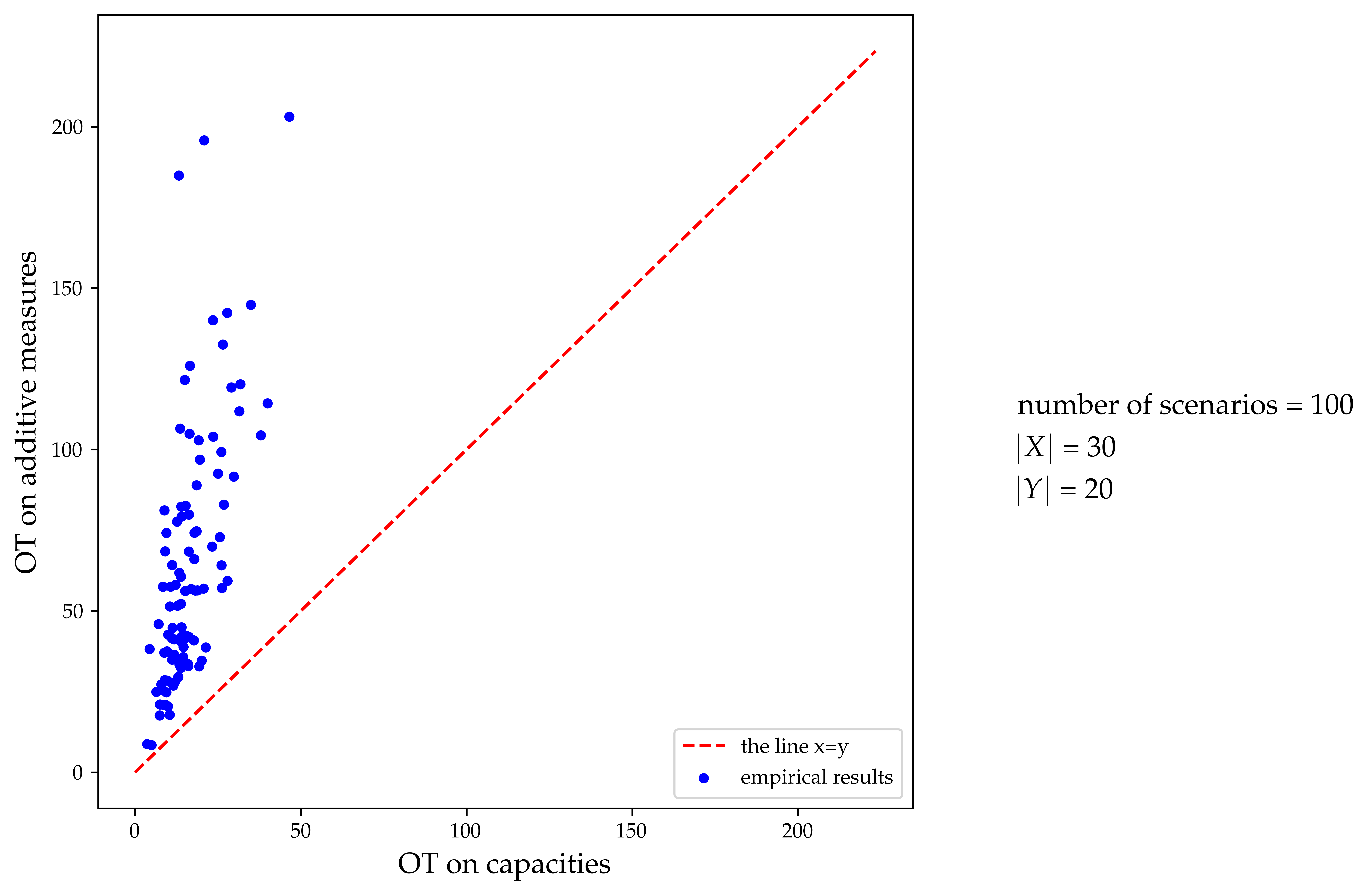}
		\caption{plot with the line $y=x$}
		\label{fig:test1a}
	\end{subfigure}
	\hfill
	\begin{subfigure}[b]{0.47\textwidth}
		\centering
		\includegraphics[width=\textwidth]{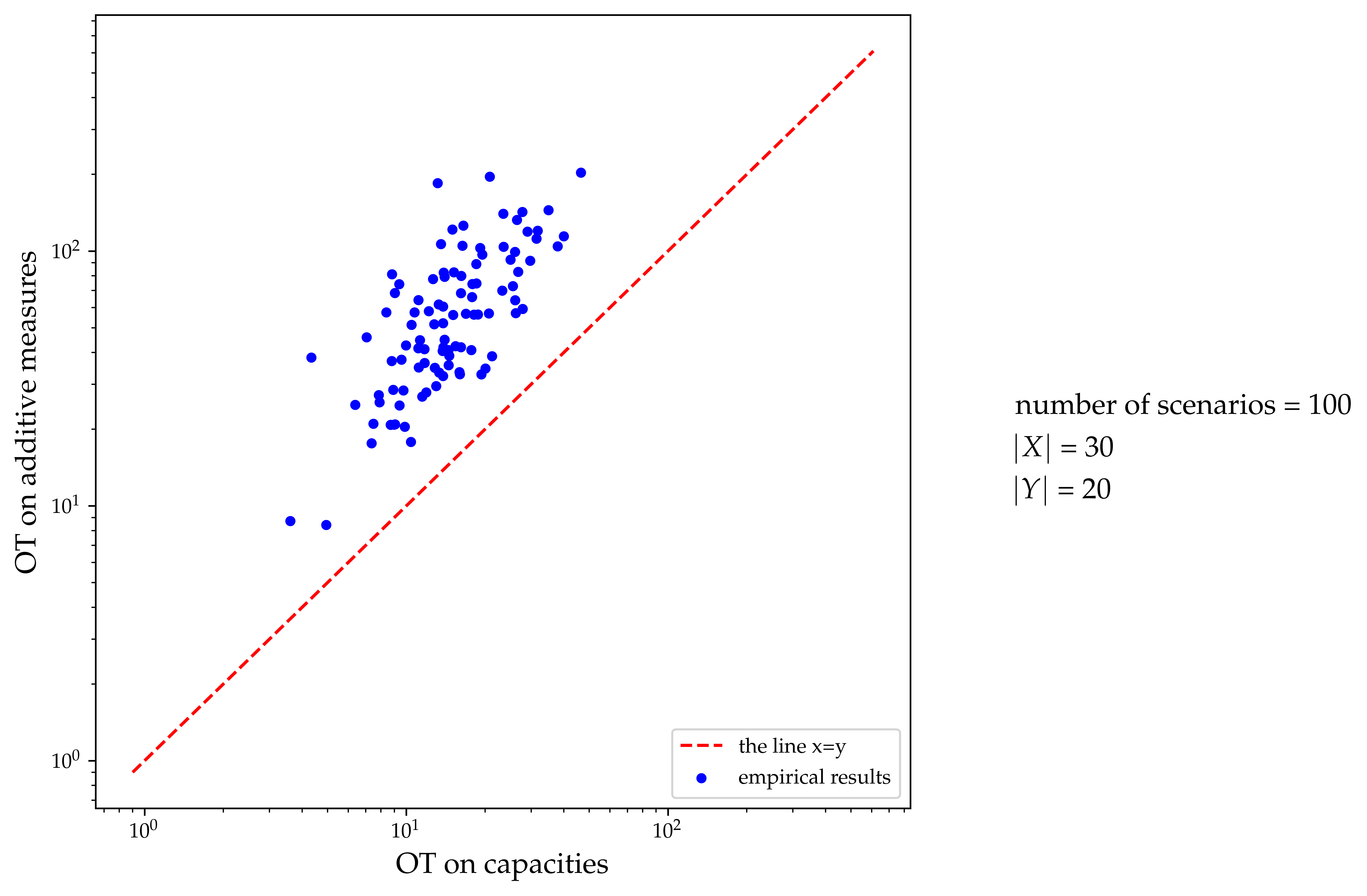}
		\caption{log-log plot with the line $y=x$}
		\label{fig:test1b}
	\end{subfigure}
	\caption{Comparison of the optimal values of the optimal transport problem for capacities with those of the optimal transport problem for additive measures, when both marginals are additive measures.}
	\label{fig:test1}
\end{figure}

% \begin{figure}[h!]
	% 	\centering
	%          \includegraphics[width=.7\textwidth]{test1_X30Y20.png}
	% 	\caption{Compare the optimal values of two minimization problems with the same marginals.}	
	
	%  \label{fig:test2} 
	% \end{figure}

\vspace{0.125cm}

For the case $|\cX|=30$ and $|\cY| = 20$, we run the following experiment 100 times. We consider the spaces $\cX = \{ 1, 2, 3, ..., 30 \}$ and $\cY = \{0, 2.2, 4.4, 6.6, ..., 41.8\}$, as well as the cost function $c(x, y) = (x-y)^2$.
{%\color{blue}
	To determine the marginal capacity $\mu$, we simulate $|\cX|-1=29$ independent random variates from a uniform distribution on $[0,1]$, and we let $U_{(i)}, i=1,\ldots,29$ be their order statistics (so that $U_{(i)}\leq U_{(i+1)}$, and $U_{(1)}$ is the smallest observation). Then set $U_{(0)}=0$  and $U_{(30)}=1$, and $\mu(\{i\}) = U_{(i)}-U_{(i-1)}$, for $i=1,\ldots,30$. The $\mu$ capacity (measure) of any other subset of $\cX$ is determined by additivity. An analogous method is used to simulate the marginal capacity  $\nu$ on $\cY$.} For each pair of simulated capacities $(\mu,\nu)$ generated in this fashion, we calculate the minima for both the optimal transport problem for measures and the one for capacities, and we compare the resulting optimal values. 
%first sample cumulative distribution functions (CDF) of each marginal distribution, which is drawn from the uniform distribution on $[0,1]$. %using numpy.random.uniform. 
%That is, we draw 29 ($|\cX|-1$) points from the uniform distribution on $[0,1]$ 
%to determine the CDF of $\mu$, then we calculate the probability density function of $\mu$ %on $\cX$ from the CDF.

\vspace{0.125cm}

The horizontal coordinates of the blue dots in Figure \ref{fig:test1a} represent the optimal values of the optimal transport problem for capacities; the vertical coordinates of the blue dots represent the minimum values of the classical optimal transport problem. 

\vspace{0.125cm}

We observe a trend that the greater the distance between the two marginal distributions, the larger the ratio between OT minimum for measures over the OT minimum for capacities. This trend is better revealed by the log-log plot in Figure \ref{fig:test1b}, showing that the one optimal value appears to behave roughly like a power of the other. The difference between these two minima implies that the classical optimal transport minimum over probability measures is inaccurate in approximating the optimal transport minimum for capacities. 

%In the following, we will use another example to elaborate on how our capacity version of the optimal transport problem is applied in counterparty credit risk.

\vspace{0.125cm}

\subsection{An Application in Counterparty Credit Risk}\label{subsec:numerical_2}

%\vspace{0.125cm}
We consider a basic model in counterparty credit risk, similar to the one  used in \cite{GHS2023}. 
%(See the next subsection for the model comparison).
%\vspace{0.125cm}
Consider a bank that trades with two counterparties whose credit exposures and the credit ratings at the end of the year determine the counterparty credit risk losses of the bank over the next year. For simplicity, we assume that there are four credit ratings, A, B, C, and D (default), with the transition probabilities in Table~\ref{tab:TransitionMatrixTable}.\footnote{Table \ref{tab:TransitionMatrixTable} is borrowed from \cite{HardySaundersQERM}.}

\begin{table}[!h]
	\begin{tabular}{|c||c|c|c|c|}\hline 
		Initial State 	& \multicolumn{4}{|c|}{Year End State} \\ \hline
		& A & B & C & D \\ \hline
		A & 0.990 & 0.007 & 0.002 & 0.001 \\ \hline
		B & 0.030 & 0.950 & 0.015 & 0.005 \\ \hline
		C & 0.015 & 0.020 & 0.960 & 0.005  \\ \hline
		D & 0 & 0 & 0 & 1\\ \hline
	\end{tabular}
	\captionsetup{justification=centering}
	\caption{Transition probabilities for a simplified credit rating system.\label{tab:TransitionMatrixTable}}
\end{table}

\vspace{0.125cm}

Assume that the initial credit ratings of counterparties 1 and 2 are $B$ and $C$, respectively. Due to ambiguity, we assume the joint rating $Y = (Y_1, Y_2)$ of these two counterparties at the year-end is represented by a capacity on $\cY$ of the form $g \circ \p$ where $g(x): = x^s$ is a concave distortion function with $s \in (0,1]$
\footnote{In particular, when $s = 1$, this capacity is the same as the additive measure $P$.}, and $\p$ is the law of  a joint probability distribution with a Gaussian copula. In particular, we let $V = (V_1, V_2)$ be a two-dimensional Gaussian random vector with mean 0 and covariance matrix
\begin{equation}
	\Sigma_{V} = \begin{pmatrix}
		1 & \rho_y\\
		\rho_y & 1
	\end{pmatrix},
	\nonumber 
\end{equation}
and define $Y_{i} = F_{i}^{\leftarrow}(\Phi(V_{i}))$, $i=1,2$, where $F_{i}$ is the marginal 
cumulative distribution function of $Y_{i}$. In particular, we have: 
\begin{equation}
	Y_1 = \begin{cases}
		D, &\text{if\ } V_1 \le \Phi^{-1}(0.005);\\
		C, &\text{if\ } \Phi^{-1}(0.005) \le V_1 \le \Phi^{-1}(0.02);\\
		B, &\text{if\ } \Phi^{-1}(0.02) \le V_1 \le \Phi^{-1}(0.97);\\
		A, &\text{if\ } V_1 \ge \Phi^{-1}(0.97),\\
	\end{cases}
\end{equation}
where $\Phi$ is the standard normal cumulative distribution function. $Y_{2}$ is 
defined similarly. 

\vspace{0.125cm}

The cardinality of $\cY$ is $16$. The probability of each pair of credit ratings can be calculated using the bivariate Gaussian distribution. For example
\begin{equation} 
	P(Y_{1}=D,Y_{2}=D) = \Phi_{2}(\Phi^{-1}(0.005),\Phi^{-1}(0.005);\rho_{y}), \nonumber
\end{equation} 
where $\Phi_{2}$ is the bivariate normal cumulative distribution function.
%We simulate $100,\!000$ scenarios from the 
%Gaussian copula in order to generate a sample of $Y_{1},Y_{2}$.
%For each subset $E \in \cY$, its associated capacity is $g(P(E))$, where $P(E)$ can be %approximated based on $100,000$ samples (i.e., $P(E)$ is approximated by the ratio of the %samples located in $E$ over the total number of samples).{\color{blue} In this example, can %we not just use the exact joint distribution? For example: 

	\vspace{0.125cm}
	
	We assume that each counterparty exposure has a (marginal) binomial distribution. In particular, we suppose that counterparty 1 has exposure $X_1$ that follows binomial($n_1$, $p_1$) and counterparty 2 has exposure $X_2$ that follows binomial($n_2$, $p_2$). The random vector $(X_{1},X_{2})$ is taken to have a Gaussian copula with correlation $\rho_{x}$. We denote the corresponding probability distribution on $\cX$ by $\Q$, and we assume that the marginal capacity $\mu=\Q$ (i.e., there is no distortion, or for the exposure capacity $s=1$). 
	
	\vspace{0.125cm}
	
	We take $n_1 = 40, p_1 = 0.4, n_2 = 25$, and $p_2 = 0.7$. % updated for both experiments!
	Then the cardinality of $\cX$ is $(n_1 + 1)\cdot (n_2+1) = 1066$. 
	%For the sake of accuracy, the sampling size for the Gaussian copula we use in this %subsection is $200,000$.
	Again, here the joint probabilities %that $X_{1}=n$ and $X_{2}=m$
	$P\left((X_1, X_2) = (n, m)\right)$ % here we use the notion of $P$ again to distinguish the capacity and probability for $X_{1}=n$ and $X_{2}=m$
	for $0\leq n \leq n_{1}$ and $0\leq m\leq n_{2}$ can be calculated using the cumulative 
	bivariate Gaussian distribution.

	\vspace{0.125cm}
	
	Finally, we describe the loss function, which is the sum of the losses due to the credit migrations of each counterparty:
	\begin{equation}
		L(X, Y) = X_1 \cdot h(Y_1) + X_2 \cdot h(Y_2), \label{LossFunctionEquation}
	\end{equation}
	where the function $h: \{A, B, C, D\} \rightarrow [0,1]$ represents the fraction of total exposure that will be lost in the next year, given the credit rating at the year-end. In this example, we take $h(A) = 0$, $h(B) = 0.1$, $h(C) = 0.2$, $h(D) = 1$ (default).
	
	\vspace{0.125cm}
	
	We look for the maximum risk represented by a Choquet integral of the loss function $L$ against a capacity $\gamma$ with prescribed marginal capacities $\mu$ and $\nu = g\circ \p$, as described above above.  That is,
	\begin{equation}\label{eqn:example2}
		\max_{\gamma \in \Pi(\mu,\nu)} \gamma(L) = \max_{\gamma\in \Pi(\mu,g\circ P)} \gamma(L).
	\end{equation}
	Note that, unlike the optimization problems in Section \ref{subsec:numerical_1}, one of the given marginals in \eqref{eqn:example2} is non-additive.

	\begin{figure}[hbt!]
		\centering
		\begin{subfigure}[b]{0.66\textwidth}
			\centering
			\includegraphics[width=\textwidth]{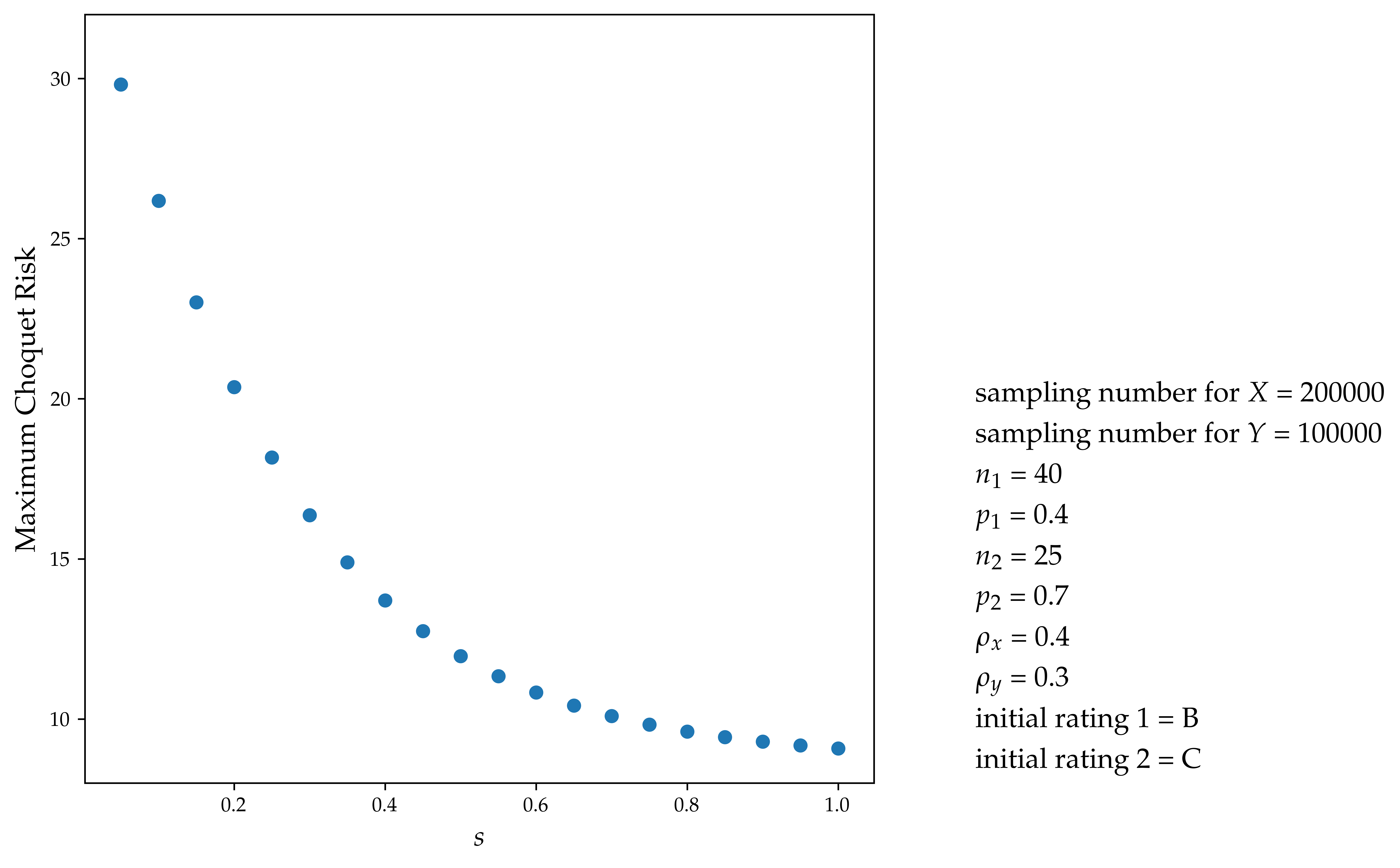}
			\caption{Maximum Choquet risk as the 
				distortion parameter $s$ varies}
			\label{fig:test2a}
		\end{subfigure}
		\hfill
		\begin{subfigure}[b]{0.47\textwidth}
			\centering
			\includegraphics[width=\textwidth]{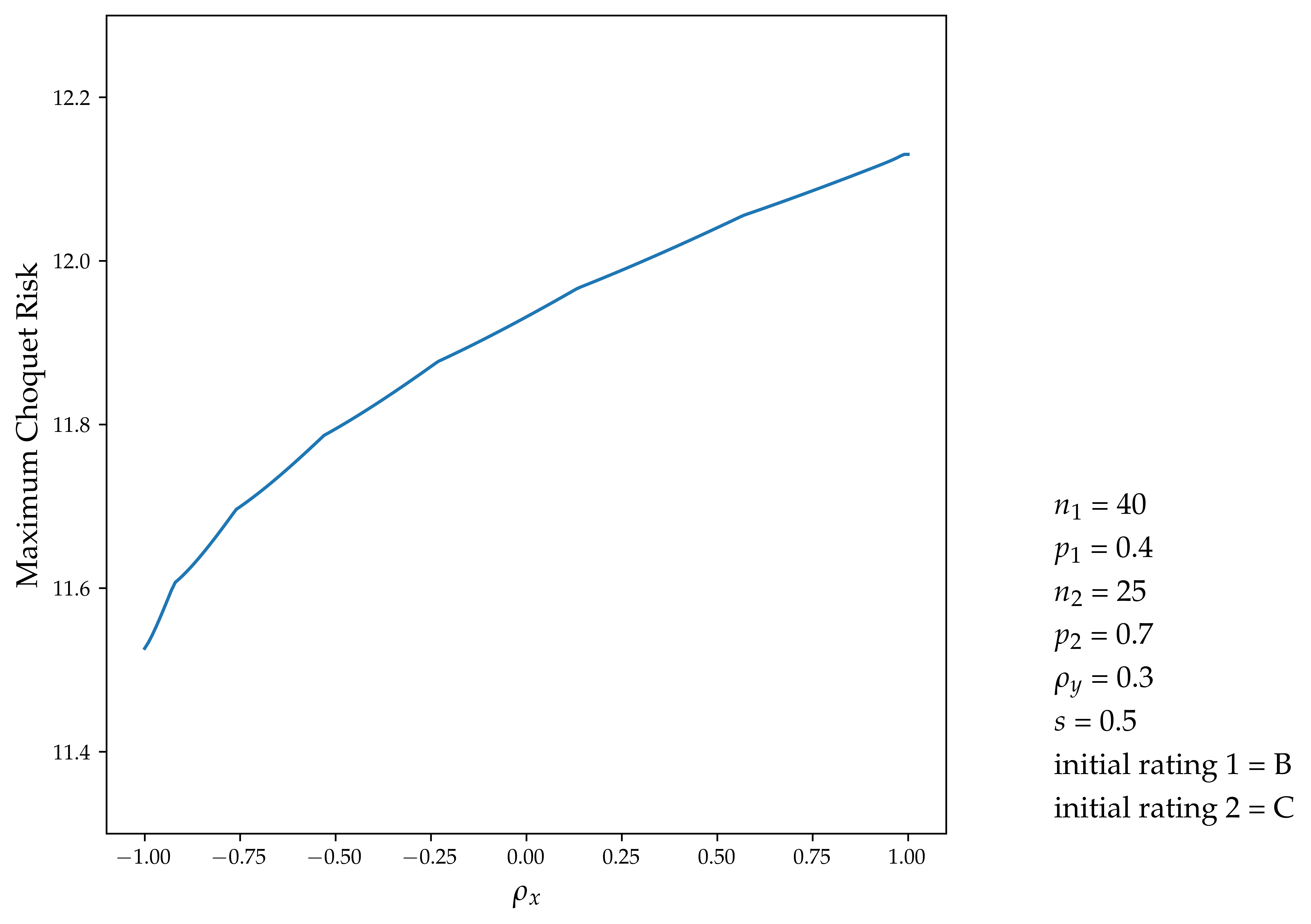}
			\caption{Maximum Choquet risk as %correlation parameter 
				$\rho_x$ varies}
			\label{fig:test2b}
		\end{subfigure}
		\hfill
		\begin{subfigure}[b]{0.47\textwidth}
			\centering
			\includegraphics[width=\textwidth]{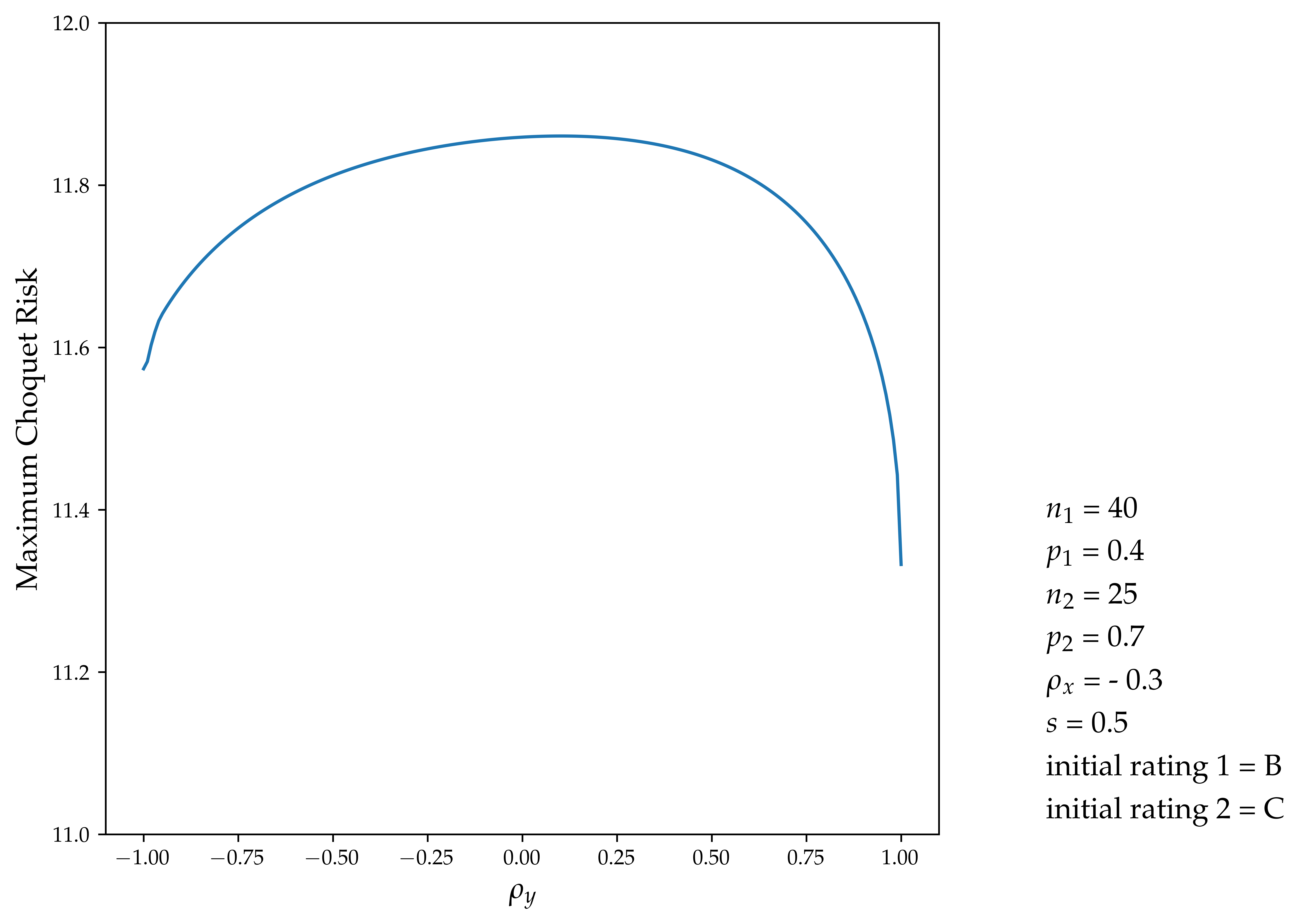}
			\caption{Maximum Choquet risk as %correlation parameter 
				$\rho_y$ varies}
			\label{fig:test2c}
		\end{subfigure}
		\caption{Sensitivity Analysis of the Maximum Choquet Risk for the 
			counterparty credit risk example.}
		\label{fig:test2}
	\end{figure}
	\vspace{0.125cm}
	
	Figure \ref{fig:test2} shows how the maximum varies along with changes in (a) the power $s$ in the distortion function, (b) the correlation $\rho_x$, or (c) the correlation $\rho_y$ in the Gaussian copula. We observe that the maximum risk is a decreasing and convex function of the distortion parameter $s$. 
	This accords with intuition, as the smaller the parameter $s$, the greater uncertainty there is regarding the marginal distribution of the credit risk factors. With the other parameters fixed, the maximum Choquet risk is an increasing and concave function of the correlation in the copula defining the distribution of the exposure factors $\mu$. Again, this makes financial sense given the nature of our loss function. If the exposures were negatively correlated, then an increase in the first term in the loss function $L$ in~(\ref{LossFunctionEquation}) would tend to correspond to a decrease in the second term. This diversification effect is amplified the greater the magnitude of the negative correlation. Similarly, if $\rho_{x}$ is large and positive, then the tail of the losses will be fatter due to the tendency for large exposures to both counterparties to occur simultaneously.

	\vspace{0.125cm}
	
	Perhaps most interesting is Figure~\ref{fig:test2c}, which considers the impact of the correlation parameter of the copula of the credit risk factors on the maximum Choquet risk, as it reveals the nature (and underlying conservatism) of the Choquet risk measure corresponding to the capacity $\pi^{*}$. The most striking aspect of the plot in Figure~\ref{fig:test2c} is that it is not monotone; there is an interior maximum of $\pi^{*}(L)$ as a function of $\rho_{y}$. To understand this, we consider explicitly a simplified version of the model. In particular, we take $X_{1}$ to be binomial with parameters $n_{1}=2$ and $p_{1}=0.4$, $X_{2}$ to be binomial with $n_{2}=2$ and $p_{2}=0.7$, and $\rho_{x}=-0.3$. Based on this specification, we can compute the probabilities for the joint distribution of the exposures $(X_{1},X_{2})$ given in Table~\ref{tab:JointExposureProbTable}.
	\begin{table}[!h]
		\begin{tabular}{|c|c|c|c|}\hline 
			& $X_{2}=0$ & $X_{2}=1$ & $X_{2}=2$ \\ \hline
			$X_{1}=0$ & 0.09 & 0.27 & 0 \\ \hline
			$X_{1}=1$ & 0 & 0.15 & 0.33 \\ \hline
			$X_{1}=2$ & 0 & 0 & 0.16  \\ \hline
		\end{tabular}
		\captionsetup{justification=centering}
		\caption{Joint exposure probabilities for the simplified counterparty credit risk example.\label{tab:JointExposureProbTable}}
	\end{table}
	
	\vspace{0.125cm}
	
	We further simplify the model by assuming only two credit states, default and no default, with both firms starting in the no-default state, and with default probability $PD=0.005$. The probability of  both counterparties defaulting together is then 
	\begin{equation}
		p_{DD}(\rho_{y}) = \Phi_{2}(\Phi^{-1}(0.005),\Phi^{-1}(0.005);\rho_{y}),
		\label{DoubleDefaultProb}
	\end{equation}
	while the probability of at least one of the counterparties defaulting is 
	$p_{D}=0.01-p_{DD}(\rho_{y})$. 
	
	\vspace{0.125cm}
	
	Let $U^{t} = \{L \geq t\}$, so that $\pi^{*}(U^{t}) = \min(\mu(U^{t}_{\cX}),\nu(U^{t}_{\cY})) 
	= \min(Q(U^{t}_{\cX}),\sqrt{P(U^{t}_{\cY})})$. Here $U^{t}_{\cX}$ is the set of $(x_{1},x_{2})$ for which there is {\em some} scenario for the credit factor $Y$ such that $L(X,Y)\geq t$. Since we can take the credit scenario to be as extreme as possible (both counterparties default), reflecting the inherent conservatism in $\pi^{*}$, we see that $U^{t}_{\cX} = \{ (x_{1},x_{2}) : x_{1} + x_{2} \geq t\}$. Simple calculations with the bivariate normal distribution with $\rho_{x}=-0.3$ then lead to the data 
	in Table~\ref{tab:ExposureProjSetsTable}.
	
	\begin{table}[!h]
		\begin{tabular}{|c|c|c|}\hline 
			$t$ values & $U^{t}_{\cX}$ & $\mu(U^{t}_{\cX})$ \\ \hline
			$t > 4$ & $\emptyset$ & 0 \\ \hline
			$3 < t \leq 4$ & $\{(2,2)\}$ & 0.0494 \\ \hline
			$2 < t \leq 3$ & $\{(1,2),(2,1),(2,2)\}$ & 0.35 \\ \hline
			$1 < t \leq 2$ & $\{ (1,0), (0,1), (0,0)\}^{c}$ & 0.8162 \\ \hline 
			$0 < t \leq 1$ & $\{(0,0)\}^{c}$ & 0.9843 \\ \hline 
			$t\leq 0$ & $\cX$ & 1 \\ \hline
		\end{tabular}
		\captionsetup{justification=centering}
		\caption{Sets $U^{t}_{\cX}$ and their capacities for the 
			simplified counterparty credit risk example.\label{tab:ExposureProjSetsTable}}
	\end{table}
	
	\vspace{0.125cm}
	
	Similarly, when considering $U^{t}_{\cY}$, we take the worst-case exposure scenario $X_{1}=X_{2}=2$, and find that $U^{t}_{\cY} = \{(y_{1},y_{2}) : h(y_{1}) + h(y_{2}) \geq \tfrac{t}{2}\}$. Recalling that $\nu(U^{t}_{\cY}) = \sqrt{P(U^{t}_{\cY})}$, we obtain the data in Table~\ref{tab:CreditProjSetsTable}.

	\begin{table}[!h]
		\begin{tabular}{|c|c|c|}\hline 
			$t$ values & $U^{t}_{\cY}$ & $\nu(U^{t}_{\cY})$ \\ \hline
			$t > 4$ & $\emptyset$ & 0 \\ \hline
			$2 < t \leq 4$ & $\{(D,D)\}$ & $\sqrt{p_{DD}(\rho_{y})}$ \\ \hline
			$0 < t \leq 2$ & $\{(A,A)\}^{c}$ & $\sqrt{0.01-p_{DD}(\rho_{y})}$ \\ \hline
			$t \leq 0$ & $\cY$ & 1 \\ \hline 
		\end{tabular}
		\captionsetup{justification=centering}
		\caption{Sets $U^{t}_{\cY}$ and their capacities for the 
			simplified counterparty credit risk example.\label{tab:CreditProjSetsTable}}
	\end{table}
	
	\vspace{0.125cm}
	
	A simple calculation then yields: 
	\begin{align} 
		\max_{\gamma\in\Pi_{\Ch}(\mu,\nu)}\gamma(L) &= 
		\int_{0}^{\infty} \pi^{*}(L\geq t) \, dt \nonumber \\
		&= \int_{0}^{4} \min\left(\mu(U^{t}_{\cX}),\nu(U^{t}_{\cY}\right)\, dt \nonumber \\
		&= 2\sqrt{0.01 - p_{DD}(\rho_{y})} + \sqrt{p_{DD}(\rho_{y})} + 0.0494, \nonumber
	\end{align} 
	and it can be seen that this function has an interior maximum (as a function of $\rho_{y}$ on $[-1,1]$). It is interesting to note that this behaviour depends on the parameters of our model, such as the probabilities of the most extreme exposure and credit scenarios. For example, with $\rho_{x}=1$ instead of $\rho_{x}=-0.3$, similar calculations give that $\pi^{*}(L)=2(\sqrt{0.01 - p_{DD}(\rho_{y})} + \sqrt{p_{DD}(\rho_{y})})$, which is monotone increasing in $\rho_{y}$.

	\vspace{0.4cm}
	
	\subsection{Comparison of Maximum Expected Shortfall and Maximum Choquet Risk with Expected Shortfall Marginal Risks} \label{subsec:numerical_3}
	%{\color{red} One can check that the functional defined above is a risk measure in terms of the random variable $L$.}
	%\\
	%\medskip
	In this subsection, we will compare the Choquet risk measure defined in the current paper with the Maximum Expected Shortfall ($\mathrm{MES}$) studied in \cite{GHS2023}.
	
	%with the maximum loss among capacities with the same given marginals via the %counterparty credit risk example described in the above section.
	
	\vspace{0.125cm}
	
	For a given loss random variable $L$ defined on $\cX\times\cY$, and for prescribed marginal probability measures $\mu$ on $\cX$ and $\nu$ on $\cY$, the {\it maximum expected shortfall} at confidence level $\alpha$ associated with $L$ is defined as
	\begin{equation*}
		\mathrm{MES}_{\alpha}(L) :=\sup_{\pi\in \Pi_a(m,n)} \mathrm{ES}_{\alpha,\pi}(L),
	\end{equation*} 
	where $\mathrm{ES}_{\alpha,\pi}$ is the expected shortfall with respect to the probability measure $\pi\in\Pi_{a}(\mu,\nu)$. In contrast to the maximum Choquet risk measure problem studied in this paper, when determining $\mathrm{MES}_{\alpha}$: 
	\vspace{0.125cm}
	\begin{itemize} 
		\item The marginal probability distributions of the risk factors on $\cX$ and $\cY$ are assumed to be known with certainty (in contrast to the case of marginal capacities, which may represent ambiguity about these marginal distributions). 
		\vspace{0.125cm}
		\item The joint risk measure is restricted to be the expected shortfall computed with respect to some probability measure $\pi\in\Pi_{a}$ (in contrast to the maximum Choquet risk measure problem, in which we consider all possible Choquet risk measures on $\cX\times\cY$ that match the given marginal Choquet risk measures on $\cX$ and $\cY$).
	\end{itemize}
	
	\vspace{0.125cm}
	
	Since expected shortfall is a distortion risk measure, the MES can be written as:
	\begin{equation}
		\mathrm{MES}_{\alpha}(L) = \sup_{\pi \in \Pi_{a}(\mu, \nu)} \int L d g_{\alpha}(\pi),
	\end{equation}
	where 
	\begin{equation}
		g_{\alpha}(x) = \begin{cases} \frac{x}{1-\alpha}, & x \in [0, 1-\alpha),\\
			1, & x \in [1-\alpha, 1],
		\end{cases} \nonumber
	\end{equation}
	is the corresponding distortion function.
	
	\vspace{0.125cm}
	
	Explicitly, the Choquet Maximum Expected Shortfall ($\mathrm{MES}$) can be defined as the maximum Choquet integral of the loss function against capacities with the same marginals as $g_{\alpha}(\pi)$.
	\begin{equation}
		\mathrm{CMES}_{\alpha}(L):= \sup_{\gamma \in \Pi(g_{\alpha}(\mu), g_{\alpha}(\nu))} \int L d\gamma.
	\end{equation}
	
	\vspace{0.125cm}
	
\begin{figure}[hbt!]
		\centering
		\begin{subfigure}[b]{0.47\textwidth}
			\centering
			\includegraphics[width=\textwidth]{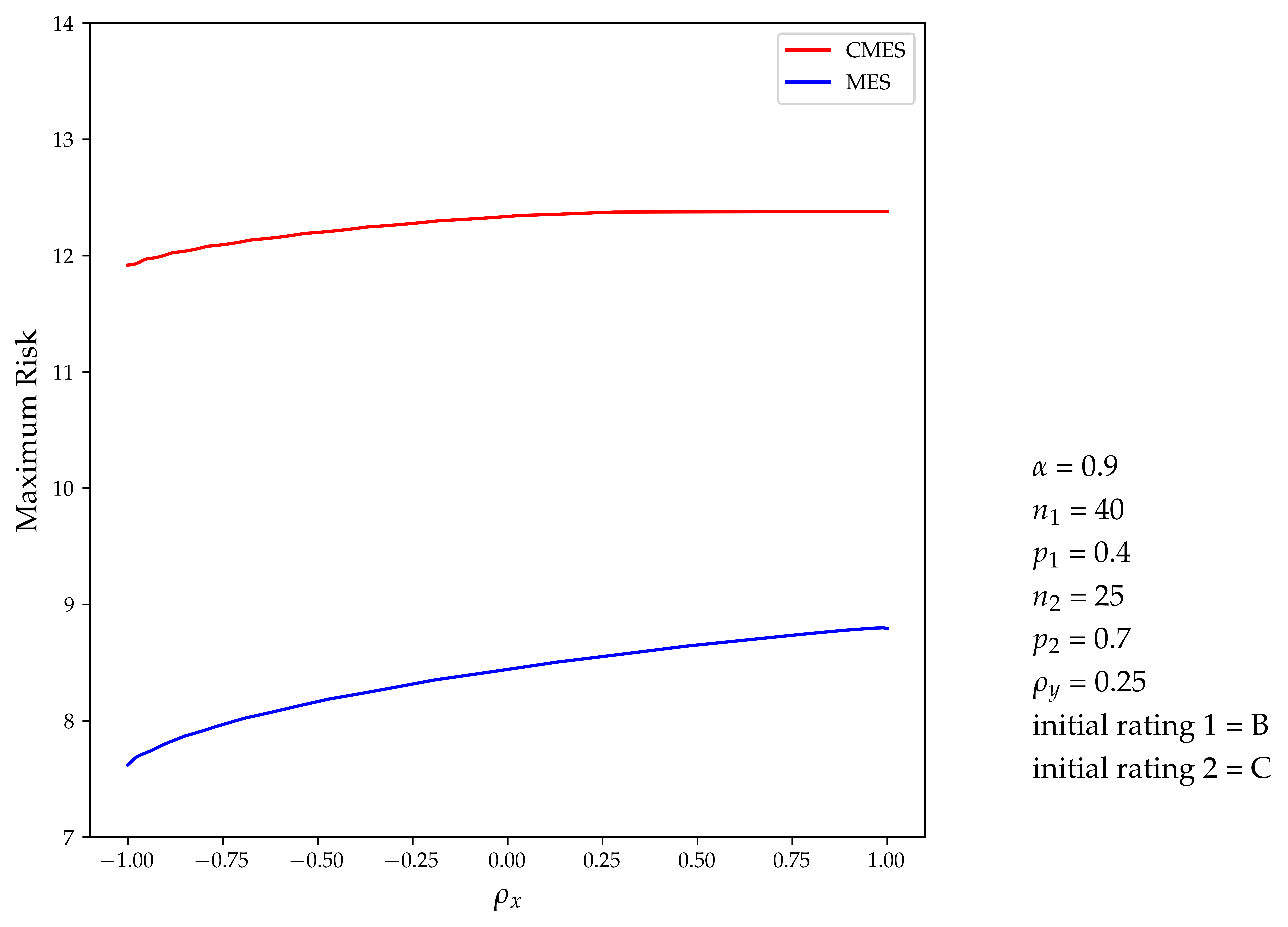}
			\caption{MES v.s. CMES with different $\rho_x$}
			\label{fig:test3a}
		\end{subfigure}
		\hfill
		\begin{subfigure}[b]{0.47\textwidth}
			\centering
			\includegraphics[width=\textwidth]{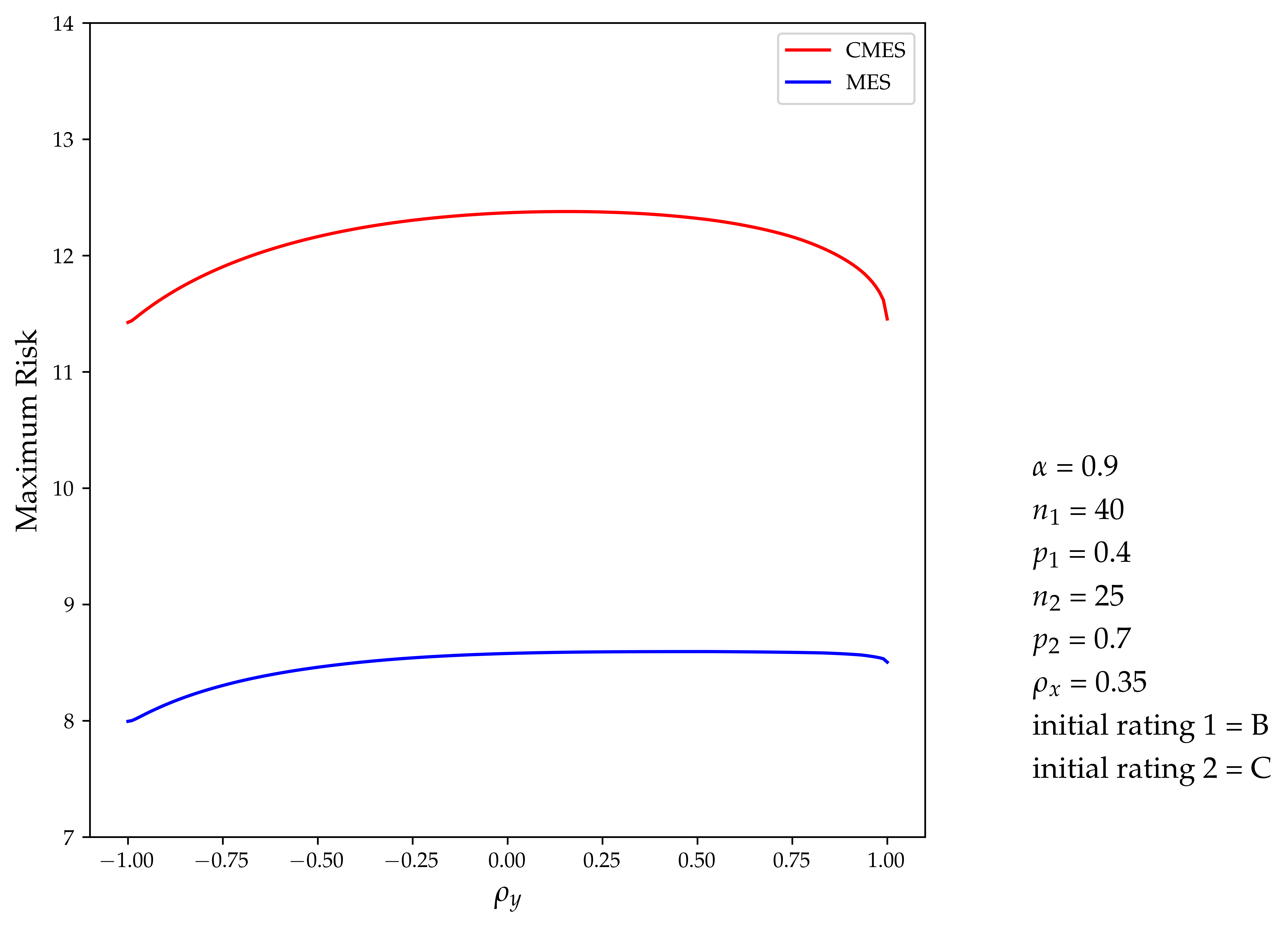}
			\caption{MES v.s. CMES with different $\rho_y$}
			\label{fig:test3b}
		\end{subfigure}
		\caption{Comparison between the
			Maximum Expected Shortfall (MES) and the
			Choquet Maximum Expected Shortfall (CMES) using different levels of the correlation parameters 
			$\rho_{x}$ and $\rho_{y}$.}
		\label{fig:test3}
	\end{figure}

	Since the feasible set for the maximum Choquet risk measure problem contains the feasible set for the maximum expected shortfall problem, we have that $\mathrm{CMES}_{\alpha}(L) \ge \mathrm{MES}_{\alpha}(L)$. In Figure \ref{fig:test3}, we compare the values of $\mathrm{CMES}_{\alpha}(L)$ and $\mathrm{MES}_{\alpha}(L)$ for {the loss random variable $L$} in the counterparty credit risk example described in the above subsection with different $\rho_x$ and $\rho_{y}$.

	\vspace{0.125cm}
	
	Throughout the experiments, we fix $\alpha = 0.9$. Assume that the counterparty 1 has initial rating B and exposure $X_1$ that follows binomial($40$, $0.4$) and that counterparty 2 has initial rating C and exposure $X_2$ that follows binomial($25$, $0.7$). In Figure \ref{fig:test3a}, the probability $\nu$ corresponds to the law of joint rating $Y = (Y_1, Y_2)$, which can be calculated using bivariate Gaussian distribution with correlation factor $\rho_y = 0.25$; similarly, the probability $\mu$ corresponds to the law of counterparty exposures $X = (X_1, X_2)$, which can be determined using bivariate Gaussian distributions with correlation factor $\rho_x$ varying from $-1$ to $1$. We plot both risk measures over different correlation factor $\rho_x$. When fixing $\rho_x = 0.35$ and allowing $\rho_y$ change from $-1$ to $1$, we obtain Figure \ref{fig:test3b}.

	\vspace{0.125cm}
	
	From the figures, one can also observe that the ratio of $\mathrm{CMES}_{0.9}(L)$ over $\mathrm{MES}_{0.9}(L)$ is between 130\% to 160\%. This ratio depends on the parameter $\alpha$ and the two given distributions, $\mu$ and $\nu$, which are eventually determined by the parameters $n_1$, $p_1$, $n_2$, $p_2$, $\rho_x$, $\rho_y$, and the values in Table \ref{tab:TransitionMatrixTable}. %For instance, if one replaces $p_1 = 0.4$ by $0.3$ in the binomial model for exposure of counterparty 1, the ratio will drop.

	\vspace{0.4cm}
	%\newpage
	%====================================================================================
	%====================================================================================
	%====================================================================================

	\section{Conclusion} \label{section:conclusion}
	
	This paper investigates the problem of bounding a Choquet risk measure of a nonlinear function of two risk factors. 
	%when the marginal distributions of the risk factors are ambiguous and represented by nonadditive measures on the marginal spaces. 
	Specifically, we assume given (marginal) capacities on the marginal spaces, representing the ambiguous distributions of the risk factors, and we consider the problem of finding the joint capacity on the product space with these given marginals, which maximizes or minimizes the Choquet integral of a given portfolio loss function. 
	
	\vspace{0.125cm}
	
	We treat this problem as a generalization of the optimal transport problem to the setting of nonadditive measures. We provide explicit characterizations of the optimal solutions for finite marginal spaces, and we investigate some of their properties. Furthermore, we investigate the relationship between properties of the marginal capacities and those of the optimizers (and, more generally, capacities in the feasible set). In particular, we show that the minimizing capacity $\pi_{*}$ is balanced if and only if both marginal capacities are balanced, and we describe its core explicitly in that case. In contrast, in all but the most trivial cases, the maximizing capacity $\pi^{*}$ is not balanced. 
	
	\vspace{0.125cm}
	
	We further discuss the connections with linear programming, showing that the optimal transport problems for capacities are linear programs, and we also characterize their duals explicitly. We investigate a series of numerical examples, including a comparison with the classical optimal transport problem, and applications to counterparty credit risk.
	
	% In the case where both marginal spaces consist of finitely many points, we showed that each of the maximization and minimization problems always has an optimal solution by providing two characterizations: one by set-wise maximization/minimization among all feasible capacities; and the other through an explicit formula involving projections of sets and the marginal capacities. 

	%====================================================================================
	%====================================================================================
	%====================================================================================
	\vspace{0.8cm}
	%\newpage
	
	\bibliographystyle{apacite}
	\bibliography{references0515.bib}
	
	\vspace{0.4cm}
	%========
\end{document}